\documentclass[11pt]{article}
 
\usepackage{amsmath, amsfonts, amssymb, amsthm,  graphicx, mathtools, enumerate}
\usepackage{mathrsfs}
\usepackage[blocks, affil-it]{authblk} 

\usepackage[numbers, square]{natbib}
\usepackage[CJKbookmarks=true,
            bookmarksnumbered=true,
			bookmarksopen=true,
			colorlinks=true,
			citecolor=red,
			linkcolor=blue,
			anchorcolor=red,
			urlcolor=blue]{hyperref} 
\usepackage[usenames]{color}

\usepackage[letterpaper, left=1.2truein, right=1.2truein, top = 1.2truein, bottom = 1.2truein]{geometry}
\usepackage[ruled, vlined, lined, commentsnumbered]{algorithm2e}

\usepackage{prettyref,soul} 
\usepackage[font=small,labelfont=it]{caption}
 \usepackage{tikz}
\usetikzlibrary{decorations.pathreplacing}
\usetikzlibrary{arrows,positioning}

\newtheorem{lemma}{Lemma}[section]

\newtheorem{thm}{Theorem}[section]

\newrefformat{eq}{(\ref{#1})}
\newrefformat{chap}{Chapter~\ref{#1}}
\newrefformat{sec}{Section~\ref{#1}}
\newrefformat{algo}{Algorithm~\ref{#1}}
\newrefformat{fig}{Fig.~\ref{#1}}
\newrefformat{tab}{Table~\ref{#1}}
\newrefformat{rmk}{Remark~\ref{#1}}
\newrefformat{clm}{Claim~\ref{#1}}
\newrefformat{def}{Definition~\ref{#1}}
\newrefformat{cor}{Corollary~\ref{#1}}
\newrefformat{lmm}{Lemma~\ref{#1}}
\newrefformat{lemma}{Lemma~\ref{#1}}
\newrefformat{prop}{Proposition~\ref{#1}}
\newrefformat{app}{Appendix~\ref{#1}}
\newrefformat{ex}{Example~\ref{#1}}
\newrefformat{exer}{Exercise~\ref{#1}}
\newrefformat{soln}{Solution~\ref{#1}}
\newrefformat{cond}{Condition~\ref{#1}}

\newcommand{\ER}{{Erd\H{o}s--R\'{e}nyi}}

\def\text#1{\mbox{\rm #1}}

\DeclarePairedDelimiter{\ceil}{\lceil}{\rceil}

\newcommand{\argmin}{\mathop{\rm argmin}}
\newcommand{\argmax}{\mathop{\rm argmax}}

\newcommand{\indc}[1]{{\mathbf{1}_{\left\{{#1}\right\}}}}
\newcommand{\norm}[1]{\|{#1} \|}

\newcommand{\wh}{\widehat}
\newcommand{\wt}{\widetilde}

\newcommand{\fnorm}[1]{\|#1\|_{\rm F}}
\newcommand{\opnorm}[1]{\|#1\|_{\rm op}}

\newcommand{\rank}{\mathop{\sf rank}}

\newcommand{\supp}{{\rm supp}}
\newcommand{\rowsupp}{{\rm rowsupp}}
\newcommand{\iprod}[2]{\left \langle #1, #2 \right\rangle}
\newcommand{\pth}[1]{\left( #1 \right)}
\newcommand{\qth}[1]{\left[ #1 \right]}

\def\X{\mathscr{X}}
\def\Z{\mathcal{Z}}
\def\T{\mathcal{T}}

\def\shape#1{
  \lower5pt\hbox{
  \hskip-7pt
  \tikzset{circ/.style={circle, draw, fill=black, scale=.2}}
  \begin{tikzpicture}[semithick,scale=.3]
  \node (l1) at (0,.866) [circ]{};
  \node (l2) at (1,.866) [circ]{};
  \node (l3) at (0.5,0) [circ]{};
  #1
  \end{tikzpicture}
  \hskip-8pt}
}

\def\edgeshape{
  \raise2pt\hbox{
  \hskip-8pt
  \tikzset{circ/.style={circle, draw, fill=black, scale=.2}}
  \begin{tikzpicture}[semithick,scale=.3]
  \node (l1) at (0,.866) [circ]{};
  \node (l2) at (1,.866) [circ]{};
  \draw[-] (l1) to node [auto] {} (l2);
  \end{tikzpicture}
  \hskip-4pt}
}

\def\veeshape{
  \lower5pt\hbox{
  \hskip-7pt
  \tikzset{circ/.style={circle, draw, fill=black, scale=.2}}
  \begin{tikzpicture}[semithick,scale=.3]
  \node (l1) at (0,.866) [circ]{};
  \node (l2) at (1,.866) [circ]{};
  \node (l3) at (0.5,0) [circ]{};
  \draw[-,color=white] (l1) to node [auto] {} (l2);
  \draw[-] (l1) to node [auto] {} (l3);
  \draw[-] (l2) to node [auto] {} (l3);
  \end{tikzpicture}
  \hskip-8pt}
}

\def\triangleshape{
  \lower5pt\hbox{
  \hskip-7pt
  \tikzset{circ/.style={circle, draw, fill=black, scale=.2}}
  \begin{tikzpicture}[semithick,scale=.3]
  \node (l1) at (0,.866) [circ]{};
  \node (l2) at (1,.866) [circ]{};
  \node (l3) at (0.5,0) [circ]{};
  \draw[-] (l1) to node [auto] {} (l3);
  \draw[-] (l2) to node [auto] {} (l3);
  \draw[-] (l1) to node [auto] {} (l2);
  \end{tikzpicture}
  \hskip-8pt}
}

\title{Minimax Rates in Network Analysis: Graphon Estimation, Community Detection and Hypothesis Testing
}
\date{~}
\author[1]{Chao Gao}
\author[2]{Zongming Ma}
\affil[1]{
University of Chicago
}
\affil[2]{
University of Pennsylvania 
}

\begin{document}
\maketitle

\begin{abstract}
This paper surveys some recent developments in fundamental limits and optimal algorithms for network analysis. 
We focus on minimax optimal rates in three fundamental problems of network analysis: graphon estimation, community detection, and hypothesis testing. 
For each problem, we review state-of-the-art results in the literature followed by general principles behind the optimal procedures that lead to minimax estimation and testing. 
This allows us to connect problems in network analysis to other statistical inference problems from a general perspective.
\smallskip

\end{abstract}


\section{Introduction}

Network analysis \citep{goldenberg2010survey} has gained considerable research interests in both theory \citep{bickel09} and applications \citep{girvan2002community,wasserman1994social}. 
In this survey, we review recent developments that establish the fundamental limits and lead to optimal algorithms in some of the most important statistical inference tasks. 
Consider a stochastic network represented by an adjacency matrix $A\in\{0,1\}^{n\times n}$. In this paper, we restrict ourselves to the setting where the network is an undirected graph without self loops. To be specific, we assume that $A_{ij}=A_{ji}\sim \text{Bernoulli}(\theta_{ij})$ for all $i<j$. The symmetric matrix $\theta\in[0,1]^{n\times n}$ models the connectivity pattern of a social network and fully characterizes the data generating process. The statistical problems we are interested is to learn structural information of the network coded in the matrix $\theta$. 
We focus on the following three problems:
\begin{enumerate}
\item \textit{Graphon estimation.} The celebrated Aldous-Hoover theorem \citep{aldous81,hoover79} asserts that the exchangeability of $\{A_{ij}\}$ implies the representation that $\theta_{ij}=f(\xi_i,\xi_j)$ with some nonparametric function $f(\cdot,\cdot)$. Here, $\xi_i$'s are i.i.d. random variables uniformly distributed in the unit interval $[0,1]$. The function $f$ is coined as the graphon of the network. The problem of graphon estimation is to estimate $f$ with the observed adjacency matrix.

\item \textit{Community detection.} Many social networks such as collaboration networks and political networks exhibit clustering structure. This means that the connectivity pattern is determined by the clustering labels of the network nodes. In general, for an assortative network, one expects that two network nodes are more likely to be connected if they are from the same cluster. For a disassortative network, the opposite pattern is expected. The task of community detection is to learn the clustering structure, and is also referred to as the problem of graph partition or network cluster analysis.

\item \textit{Hypothesis testing.} Perhaps the most fundamental question for network analysis is whether a network has some structure. 
For example, an \ER~graph has a constant connectivity probability for all edges, and is regarded to have no interesting structure. In comparison, a stochastic block model has a clustering structure that governs the connectivity pattern. Therefore, before conducting any specific network analysis, one should first test whether a network has some structure or not. The test between an \ER~graph and a stochastic block model is one of the simplest examples.
\end{enumerate}

This survey will emphasize the developments of the minimax rates of the problems. 
The state-of-the-art of the three problems listed above will be reviewed in Section \ref{sec:graphon}, Section \ref{sec:community-detection}, and Section \ref{sec:test}, respectively. 
In each section, we will introduce critical mathematical techniques that we use to derive optimal solutions.
When appropriate, we will also discuss the general principles behind the problems.
This allows us to connect the results of the network analysis to some other interesting statistical inference problems. 
%

Real social networks are often sparse, which means that the number of edges are of a smaller order compared with the number of nodes squared. How to model sparse networks is a longstanding topic full of debate \citep{lloyd2012random,bickel09,crane2016edge,caron2017sparse}. In this paper, we adopt the notion of network sparsity $\max_{1\leq i<j\leq n}\theta_{ij}=o(1)$, which is proposed by \cite{bickel09}. 
Theoretical foundations of this sparsity notion were investigated by \cite{borgs2014p1,borgs2014p2}. 
There are other, perhaps more natural, notions of network sparsity, and we will discuss potential open problems in Section \ref{sec:disc}.

We close this section by introducing some notation that will be used in the paper. For an integer $d$, we use $[d]$ to denote the set $\{1,2,...,d\}$. Given two numbers $a,b\in\mathbb{R}$, we use $a\vee b=\max(a,b)$ and $a\wedge b=\min(a,b)$. For two positive sequences $\{a_n\},\{b_n\}$, $a_n\lesssim b_n$ means $a_n\leq Cb_n$ for some constant $C>0$ independent of $n$, and $a_n\asymp b_n$ means $a_n\lesssim b_n$ and $b_n\lesssim a_n$. We write $a_n\ll b_n$ if $a_n/b_n\rightarrow 0$. For a set $S$, we use $\indc{S}$ to denote its indicator function and $|S|$ to denote its cardinality. For a vector $v\in\mathbb{R}^d$, its norms are defined by $\norm{v}_1=\sum_{i=1}^n|v_i|$, $\norm{v}^2=\sum_{i=1}^nv_i^2$ and $\norm{v}_{\infty}=\max_{1\leq i\leq n}|v_i|$. For two matrices $A,B\in\mathbb{R}^{d_1\times d_2}$, their trace inner product is defined as $\iprod{A}{B}=\sum_{i=1}^{d_1}\sum_{j=1}^{d_2}A_{ij}B_{ij}$. The Frobenius norm and the operator norm of $A$ are defined by $\fnorm{A}=\sqrt{\iprod{A}{A}}$ and $\opnorm{A}=s_{\max}(A)$, where $s_{\max}(\cdot)$ denotes the largest singular value.

\section{Graphon estimation}\label{sec:graphon}

\subsection{Problem settings}

Graphon is a nonparametric object that determines the data generating process of a random network. 
The concept is from the literature of exchangeable arrays \citep{aldous81,hoover79,kallenberg89} and graph limits \citep{lovasz12,diaconis07}. 
We consider a random graph with adjacency matrix $\{A_{ij}\}\in\{0,1\}^{n\times n}$, whose sampling procedure is determined by
\begin{equation}
(\xi_1,...,\xi_n)\sim\mathbb{P}_{\xi},\quad A_{ij}|(\xi_i,\xi_j)\sim\text{Bernoulli}(\theta_{ij}),\quad\text{where }\theta_{ij}=f(\xi_i,\xi_j).\label{eq:graphon}
\end{equation}
For $i\in[n]$, $A_{ii}=\theta_{ii}=0$. 
Conditioning on $(\xi_1,...,\xi_n)$, the $A_{ij}$'s are mutually independent across all $i<j$. 
The function $f$ on $[0,1]^2$, which is assumed to be symmetric, is called graphon. The graphon offers a flexible nonparametric way of modeling stochastic networks. We note that exchangeability leads to independent random variables $(\xi_1,...,\xi_n)$ sampled from $\text{Uniform}[0,1]$, but for the purpose of estimating $f$, we do not require this assumption.

We point out an interesting connection between graphon estimation and nonparametric regression. In the formulation of (\ref{eq:graphon}), suppose we observe both the adjacency matrix $\{A_{ij}\}$ and the latent variables $\{(\xi_i,\xi_j)\}$, then $f$ can simply be regarded as a regression function that maps $(\xi_i,\xi_j)$ to the mean of $A_{ij}$. However, in the setting of network analysis, we only observe the adjacency matrix $\{A_{ij}\}$. The latent variables are usually used to model latent features of the network nodes \citep{hoff2002latent,ma2017exploration}, and are not always available in practice. Therefore, graphon estimation is essentially a nonparametric regression problem without observing the covariates, which leads to a new phenomenon in the minimax rate that we will present below.

In the literature, various estimators have been proposed. For example, a singular value threshold method is analyzed by \cite{chatterjee12}, later improved by \cite{xu2017rates}. The paper \cite{lloyd2012random} considers a Bayesian nonparametric approach. Another popular procedure is to estimate the graphon via histogram or stochastic block model approximation \citep{wolfe13,chan14,airoldi13,olhede13,borgs2015private,borgs2015consistent}. Minimax rates of graphon estimation are investigated by \cite{gao2015rate,gao2016optimal,klopp2017structured}.

\subsection{Optimal rates}

Before discussing the minimax rate of estimating a nonparametric graphon, we first consider graphons that are block-wise constant functions. This is equivalently recognized as stochastic block models (SBMs) \citep{holland83,nowicki01}. Consider $A_{ij}\sim \text{Bernoulli}(\theta_{ij})$ for all $1\leq i<j\leq n$. The class of SBMs with $k$ clusters is defined as
\begin{eqnarray}
\label{eq:SBM-k} \Theta_k &=& \Bigg\{\{\theta_{ij}\}\in [0,1]^{n\times n}: \theta_{ii}=0,\text{ }\theta_{ij}=B_{uv}=B_{vu}\\
\nonumber && \text{ for }(i,j)\in z^{-1}(u)\times z^{-1}(v)\text{ with some }B_{uv}\in[0,1]\text{ and }z\in[k]^n\Bigg\}.
\end{eqnarray}
In other words, the network nodes are divided into $k$ clusters that are determined by the cluster labels $z$. The subsets $\{\mathcal{C}_u(z)\}_{z\in[k]}$ with $\mathcal{C}_u(z)=\{i\in[n]:z(i)=u\}$ form a partition of $[n]$. The mean matrix $\theta\in[0,1]^{n\times n}$ is a piecewise constant with respect to the blocks $\{\mathcal{C}_u(z)\times \mathcal{C}_v(z):u,v\in [k] \}$.

In this setting, graphon estimation is the same as estimating the mean matrix $\theta$. If we know the clustering labels $z$, then we can simply calculate the sample averages of $\{A_{ij}\}$ in each block $\mathcal{C}_u(z)\times \mathcal{C}_v(z)$. Without the knowledge of $z$, a least-squares estimator proposed by \cite{gao2015rate} is
\begin{equation}
\wh{\theta}=\argmin_{\theta\in\Theta_k}\fnorm{A-\theta}^2,\label{eq:least-squares}
\end{equation}
which can be understood as the sample averages of $\{A_{ij}\}$ over the estimated blocks $\{\mathcal{C}_u(\wh{z})\times \mathcal{C}_v(\wh{z})\}$.

To study the performance of the least-squares estimator $\wh{\theta}$, we need to introduce some additional notation. 
Since $\wh{\theta}\in\Theta_k$, the estimator can be written as $\wh{\theta}_{ij}=\wh{B}_{\wh{z}(i)\wh{z}(j)}$ for some $\wh{B}\in[0,1]^{k\times k}$ and some $\wh{z}\in[k]^n$. The true matrix that generates $A$ is denoted by $\theta^*$. 
Then, we define
$$\wt{\theta}=\argmin_{\theta\in\Theta_k(\wh{z})}\fnorm{\theta^*-\theta}^2.$$
Here, the class $\Theta_k(\wh{z})\subset\Theta_k$ consists of all SBMs with clustering structures determined by $\wh{z}$. Then, we immediately have the Pythagorean identity
\begin{equation}
\fnorm{\wh{\theta}-\theta^*}^2 = \fnorm{\wh{\theta}-\wt{\theta}}^2 + \fnorm{\wt{\theta}-\theta^*}^2.\label{eq:Pythagorean}
\end{equation}
By the definition of $\wh{\theta}$, we have the basic inequality $\fnorm{\wh{\theta}-A}^2\leq \fnorm{\theta^*-A}^2$. After a simple rearrangement, we have
\begin{eqnarray*}
\fnorm{\wh{\theta}-\theta^*}^2 &\leq& 2\left|\iprod{\wh{\theta}-\theta^*}{A-\theta^*}\right| \\
&\leq& 2\fnorm{\wh{\theta}-\wt{\theta}}\left|\iprod{\frac{\wh{\theta}-\wt{\theta}}{\fnorm{\wh{\theta}-\wt{\theta}}}}{A-\theta^*}\right| + 2\fnorm{\wt{\theta}-\theta^*}\left|\iprod{\frac{\wt{\theta}-\theta^*}{\fnorm{\wt{\theta}-\theta^*}}}{A-\theta^*}\right| \\
&\leq& \fnorm{\wh{\theta}-\theta^*}\sqrt{\left|\iprod{\frac{\wh{\theta}-\wt{\theta}}{\fnorm{\wh{\theta}-\wt{\theta}}}}{A-\theta^*}\right|^2+\left|\iprod{\frac{\wt{\theta}-\theta^*}{\fnorm{\wt{\theta}-\theta^*}}}{A-\theta^*}\right|^2},
\end{eqnarray*}
where the last inequality is by Cauchy-Schwarz and (\ref{eq:Pythagorean}). Therefore, we have
\begin{eqnarray*}
\fnorm{\wh{\theta}-\theta^*}^2 &\leq& \left|\iprod{\frac{\wh{\theta}-\wt{\theta}}{\fnorm{\wh{\theta}-\wt{\theta}}}}{A-\theta^*}\right|^2+\left|\iprod{\frac{\wt{\theta}-\theta^*}{\fnorm{\wt{\theta}-\theta^*}}}{A-\theta^*}\right|^2 \\
&\leq& \sup_{\{v\in\Theta_k:\fnorm{v}=1\}}|\iprod{v}{A-\theta^*}|^2 + \max_{1\leq j\leq k^n}|\iprod{v_j}{A-\theta^*}|^2,
\end{eqnarray*}
where $\{v_j\}_{1\leq j\leq k^n}$ are $k^n$ fixed matrices with Frobenius norm $1$.
To understand the last inequality above, observe that $\frac{\wh{\theta}-\wt{\theta}}{\fnorm{\wh{\theta}-\wt{\theta}}}$ belongs to $\Theta_k$ and has Frobenius norm $1$, and the matrix $\frac{\wt{\theta}-\theta^*}{\fnorm{\wt{\theta}-\theta^*}}$ takes at most $k^n$ different values. Finally, an empirical process argument and a union bound leads to the inequalities
\begin{eqnarray*}
\mathbb{E}\left[\sup_{\{v\in\Theta_k:\fnorm{v}=1\}}|\iprod{v}{A-\theta^*}|^2\right] &\lesssim& k^2+n\log k, \\
\mathbb{E}\left[\max_{1\leq j\leq k^n}|\iprod{v_j}{A-\theta^*}|^2\right] &\lesssim& n\log k,
\end{eqnarray*}
which then implies the bound
\begin{equation}
\mathbb{E}\fnorm{\wh{\theta}-\theta^*}^2 \lesssim k^2 + n\log k. \label{eq:SBM-upper}
\end{equation}

The upper bound (\ref{eq:SBM-upper}) consists of two terms. The first term $k^2$ corresponds to the number of parameters we need to estimate in an SBM with $k$ clusters. 
The second term results from not knowing the exact clustering structure. Since there are in total $k^n$ possible clustering configurations, the complexity $\log(k^n)=n\log k$ enters the error bound. Even though the bound (\ref{eq:SBM-upper}) is achieved by an estimator that knows the value of $k$, a penalized version of the least-squares estimator with the penalty $\lambda(k^2 + n\log k)$ can achieve the same bound (\ref{eq:SBM-upper}) without the knowledge of $k$.

The paper \cite{gao2015rate} also shows that the upper bound (\ref{eq:SBM-upper}) is sharp by proving a matching minimax lower bound. While it is easy to see that the first term $k^2$ cannot be avoided by a classical lower bound argument of parametric estimation, the necessity of the second term $n\log k$ requires a very delicate lower bound construction. 
It was proved by \cite{gao2015rate} that it is possible to construct a $B\in[0,1]^{k\times k}$, such that the set $\{B_{z(i)z(j)}:z\in[k]^n\}$ has a packing number bounded below by $e^{cn\log k}$ with respect to the norm $\fnorm{\cdot}$ and the radius at the order of $\sqrt{n\log k}$. This fact, together with a standard Fano inequality argument, leads to the desired minimax lower bound.

We summarize the above discussion into the following theorem.
\begin{thm}[Gao, Lu and Zhou \citep{gao2015rate}]\label{thm:sbm-minimax}
For the loss function $L(\wh{\theta},\theta)={{n\choose 2}}^{-1}\sum_{1\leq i<j\leq n}(\wh{\theta}_{ij}-\theta_{ij})^2$, we have
$$\inf_{\wh{\theta}}\sup_{\theta\in\Theta_k}\mathbb{E}L(\wh{\theta},\theta)\asymp \frac{k^2}{n^2} + \frac{\log k}{n},$$
for all $1\leq k\leq n$.
\end{thm}

Having understood minimax rates of estimating mean matrices of SBMs, 
we are ready to discuss minimax rates of estimating general nonparametric graphons.
We consider the following loss function that is widely used in the literature of nonparametric regression,
$$L(\wh{f},f)=\frac{1}{{n\choose 2}}\sum_{1\leq i<j\leq n}\left(\wh{f}(\xi_i,\xi_j)-f(\xi_i,\xi_j)\right)^2.$$
Note that $L(\wh{f},f)=L(\wh{\theta},\theta)$ if we let $\wh{\theta}_{ij}=\wh{f}(\xi_i,\xi_j)$ and $\theta_{ij}=f(\xi_i,\xi_j)$.
Then, the minimax risk is defined as
$$\inf_{\wh{f}}\sup_{f\in\mathcal{H}_{\alpha}(M)}\sup_{\mathbb{P}_{\xi}}\mathbb{E}L(\wh{f},f).$$
Here, the supreme is over both the function class $\mathcal{H}_{\alpha}(M)$ and the distribution $\mathbb{P}_{\xi}$ that the latent variables $(\xi_1,...,\xi_n)$ are sampled from. While $\mathbb{P}_{\xi}$ is allowed to range from the class of all distributions, the H\"{o}lder class $\mathcal{H}_{\alpha}(M)$ is defined as
$$\mathcal{H}_{\alpha}(M)=\left\{\|f\|_{\mathcal{H}_{\alpha}}\leq M: f(x,y)=f(y,x)\text{ for }x\geq y\right\},$$
where $\alpha>0$ is the smoothness parameter and $M>0$ is the size of the class. Both are assumed to be constants. In the above definition, $\|f\|_{\mathcal{H}_{\alpha}}$ is the H\"{o}lder norm of the function $f$ (see \cite{gao2015rate} for the details).

The following theorem gives the minimax rate of the problem.
\begin{thm}[Gao, Lu and Zhou \citep{gao2015rate}]\label{thm:graphon-minimax}
We have
$$\inf_{\wh{f}}\sup_{f\in\mathcal{H}_{\alpha}(M)}\sup_{\mathbb{P}_{\xi}}\mathbb{E}L(\wh{f},f) \asymp \begin{cases}
n^{-\frac{2\alpha}{\alpha+1}}, & 0<\alpha<1,\\
\frac{\log n}{n}, & \alpha\geq 1,
\end{cases}$$
where the expectation is jointly over $\{A_{ij}\}$ and $\{\xi_i\}$.
\end{thm}
The minimax rate in Theorem \ref{thm:graphon-minimax} exhibits different behaviors in the two regimes depending on whether $\alpha\geq 1$ or not. For $\alpha\in(0,1)$, we obtain the classical minimax rate for nonparametric regression. To see this, one can related the graphon estimation problem to a two-dimensional nonparametric regression problem with sample size $N=\frac{n(n-1)}{2}$, and then it is easy to see that $N^{-\frac{2\alpha}{2\alpha+d}}\asymp n^{-\frac{2\alpha}{\alpha+1}}$ for $d=2$. This means for a nonparametric graphon that is not so smooth, whether or not the latent variables $\{(\xi_i,\xi_j)\}$ are observed does not affect the minimax rate. In contrast, when $\alpha\geq 1$, the minimax rate scales as $\frac{\log n}{n}$, which does not depend on the value of $\alpha$ anymore. In this regime, there is a significant difference between the graphon estimation problem and the regression problem.

Both the upper and lower bounds in Theorem \ref{thm:graphon-minimax} can be derived by an SBM approximation. The minimax rate given by Theorem \ref{thm:graphon-minimax} can be equivalently written as
$$\min_{1\leq k\leq n}\left\{\frac{k^2}{n^2}+\frac{\log k}{n}+k^{-2(\alpha\wedge 1)}\right\},$$
where $\frac{k^2}{n^2}+\frac{\log k}{n}$ is the optimal rate of estimating a $k$-cluster SBM in Theorem \ref{thm:sbm-minimax}, and $k^{-2(\alpha\wedge 1)}$ is the approximation error for an $\alpha$-smooth graphon by a $k$-cluster SBM. As a consequence, the least-squares estimator (\ref{eq:least-squares}) is rate-optimal with $k\asymp n^{\frac{1}{\alpha\wedge 1+1}}$. The result justifies the strategies of estimating a nonparametric graphon by network histograms in the literature \citep{wolfe13,chan14,airoldi13,olhede13}.
 
Despite its rate-optimality, an disadvantage of the least-squares estimator (\ref{eq:least-squares}) is its computational intractability. 
A naive algorithm requires an exhaustive search over all $k^n$ possible clustering structures. Although a two-way $k$-means algorithm in \cite{gao2016optimal} works well in practice, there is no theoretical guarantee that the algorithm can find the global optimum in polynomial time. An alternative strategy is to relax the constraint in the least-squares optimization. For instance, let $\wt{\Theta}_k$ be the set of all symmetric matrices $\theta\in[0,1]^{n\times n}$ that have at most $k$ ranks. It is easy to see $\Theta_k\subset\wt{\Theta}_k$. Moreover, the relaxed estimator $\wh{\theta}=\argmin_{\theta\in\wt{\Theta}_k}\fnorm{A-\theta}^2$ can be computed efficiently through a simple eigenvalue decomposition. This is closely related to the procedures discussed in \cite{chatterjee12}. However, such an estimator can only achieve the rate $\frac{k}{n}$, which can be much slower than the minimax rate $\frac{k^2}{n^2}+\frac{\log k}{n}$. To the best of our knowledge, $\frac{k}{n}$ is the best known rate that can be achieved by a polynomial-time algorithm so far. We refer the readers to \cite{xu2017rates} for more details on this topic.

\subsection{Extensions to sparse networks}

In many practical situations, sparse networks are more useful. 
A network is sparse if the maximum probability of $\{A_{ij}=1\}$ tends to zero as $n$ tends to infinity. 
A sparse graphon $f$ is a symmetric nonnegative function on $[0,1]$ that satisfies $\sup_{x,y}f(x,y)\leq \rho=o(1)$ \citep{bickel09,borgs2014p1,borgs2014p2}. Analogously, a sparse SBM is characterized by the space $\Theta_k(\rho)=\{\theta\in\Theta_k:\max_{ij}\theta_{ij}\leq \rho\}$. An extension of Theorem \ref{thm:sbm-minimax} is given by the following result.
\begin{thm}[Klopp, Tsybakov and Verzelen \citep{klopp2017oracle}; Gao, Ma, Lu and Zhou \citep{gao2016optimal}]\label{thm:minimax-sparse-graphon}
We have
$$\inf_{\wh{\theta}}\sup_{\theta\in\Theta_k}\mathbb{E}L(\wh{\theta},\theta)\asymp \min\left\{\rho\left(\frac{k^2}{n^2} + \frac{\log k}{n}\right),\rho^2\right\},$$
for all $1\leq k\leq n$.
\end{thm}
Theorem \ref{thm:minimax-sparse-graphon} recovers the minimax rate of Theorem \ref{thm:sbm-minimax} if we set $\rho\asymp 1$. The result was obtained independently by \cite{klopp2017oracle} and \cite{gao2016optimal} around the same time. Besides the the loss function $L(\cdot,\cdot)$ on the probability matrix, the paper \cite{klopp2017oracle} also considered integrated loss for the graphon function. 

To achieve the minimax rate, one can consider the least-squares estimator
\begin{equation}
\wh{\theta}=\argmin_{\theta\in\Theta_k(\rho)}\fnorm{A-\theta}^2\label{eq:sparse-least-squares}
\end{equation}
when $\frac{k^2}{n^2} + \frac{\log k}{n}\geq \rho$. In the situation when $\frac{k^2}{n^2} + \frac{\log k}{n}< \rho$, the minimax rate is $\rho^2$ and can be trivially achieved by $\wh{\theta}=0$.

Theorem \ref{thm:minimax-sparse-graphon} also leads to optimal rates of nonparametric sparse graphon estimation in a H\"{o}lder space \citep{klopp2017oracle,gao2016optimal}. 
In addition, sparse graphon estimation in a privacy-aware setting \citep{borgs2015private} and a heavy-tailed setting \citep{borgs2015consistent} have also been considered in the literature.

\subsection{Biclustering and related problems}

SBM can be understood as a special case of biclustering. 
A matrix has a biclustering structure if it is block-wise constant with respect to both row and column clustering structures. 
The biclustering model was first proposed by \cite{hartigan1972direct}, and has been widely used in modern gene expression data analysis \citep{cheng2000biclustering,madeira2004biclustering}. 
Mathematically, we consider the following parameter space
\begin{eqnarray*}
\Theta_{k,l} &=& \Bigg\{\{\theta_{ij}\}\in \mathbb{R}^{n\times m}: \theta_{ij}=B_{z_1(i)z_2(j)}, B\in\mathbb{R}^{k\times l}, z_1\in[k]^n, z_2\in[l]^m\Bigg\}.
\end{eqnarray*}
Then, for the loss function $L(\wh{\theta},\theta)=\frac{1}{nm}\sum_{i=1}^n\sum_{j=1}^m(\wh{\theta}_{ij}-\theta_{ij})^2$, 
it has been shown in \cite{gao2015rate,gao2016optimal} that
\begin{equation}
\inf_{\wh{\theta}}\sup_{\theta\in\Theta_{k,l}}\mathbb{E}L(\wh{\theta},\theta)\asymp \frac{kl}{mn} + \frac{\log k}{m} + \frac{\log l}{n}, \label{eq:biclustering-minimax}
\end{equation}
as long as $\log k\asymp \log l$. The minimax rate (\ref{eq:biclustering-minimax}) holds under both Bernoulli and Gaussian observations. When $k=l$ and $m=n$, the result (\ref{eq:biclustering-minimax}) recovers Theorem \ref{thm:sbm-minimax}.

The minimax rate (\ref{eq:biclustering-minimax}) reveals a very important principle of sample complexity. In fact, for a large collection of popular problems in high-dimensional statistics, the minimax rate is often in the form of
\begin{equation}
\frac{(\#\text{parameters}) + \log(\#\text{models})}{\text{\#samples}}.\label{eq:general}
\end{equation}
For the biclustering problem, $nm$ is the sample size and $kl$ is the number of parameters. Since the number of biclustering structures is $k^nl^m$, the formula (\ref{eq:general}) gives (\ref{eq:biclustering-minimax}).

To understand the general principle (\ref{eq:general}), we need to discuss the structured linear model introduced by \cite{gao2015general}. In the framework of structured linear models, the data can be written as
$$Y=\X_Z(B)+W\in\mathbb{R}^N,$$
where $\X_Z(B)$ is the signal to be recovered and $W$ is a mean-zero noise. The signal part $\X_Z(B)$ consists of a linear operator $\X_Z(\cdot)$ indexed by the model/structure $Z$ and parameters that are organized as $B$. The structure $Z$ is in some discrete space $\Z_{\tau}$, which is further indexed by $\tau\in\T$ for some finite set $\T$. We introduce a function $\ell(\Z_{\tau})$ that determines the dimension of $B$. In other words, we have $B\in\mathbb{R}^{\ell(\Z_{\tau})}$. Then, the optimal rate that recovers the signal $\theta=\X_Z(B)$ with respect to the loss function $L(\wh{\theta},\theta)=\frac{1}{N}\sum_{i=1}^N(\wh{\theta}_i-\theta_i)^2$ is given by
\begin{equation}
\frac{\ell(\Z_{\tau})+\log|\Z_{\tau}|}{N}.\label{eq:general-rigorous}
\end{equation}
We note that (\ref{eq:general-rigorous}) is a mathematically rigorous version of (\ref{eq:general}). In \cite{gao2015general}, a Bayesian nonparametric procedure was proposed to achieve the rate (\ref{eq:general-rigorous}). 
Minimax lower bounds in the form of (\ref{eq:general-rigorous}) have been investigated by \cite{klopp2017structured} under a slightly different framework. 
Below we present a few important examples of the structured linear models.

\textit{Biclustering.} In this model, it is convenient to organize $\X_Z(B)$ as a matrix in $\mathbb{R}^{n\times m}$ and then $N=nm$. The linear operator $\X_Z(\cdot)$ is determined by $[\X_Z(B)]_{ij}=B_{z_1(i)z_2(j)}$ with $Z=(z_1,z_2)$. With the relations $\tau=(k,l)$, $\T=[n]\times[m]$, $\Z_{k,l}=[k]^n\times[l]^m$, we get $\ell(Z_{k,l})=kl$ and $\log|\Z_{k,l}|=n\log k+m\log l$, and the rate (\ref{eq:biclustering-minimax}) can be derived from (\ref{eq:general-rigorous}).

\textit{Sparse linear regression.} The linear model $X\beta$ with a sparse $\beta\in\mathbb{R}^p$ can also be written as $\X_Z(B)$. To do this, note that a sparse $\beta$ implies a representation $\beta^T=(\beta_S^T,0_{S^c}^T)$ for some subset $S\subset[p]$. Then, $X\beta=X_{*S}\beta_S=\X_Z(B)$, with the relations $Z=S$, $\tau=s$, $\T=[p]$, $\Z_s=\{S\subset[p]:|S|=s\}$, $\ell(\Z_s)=s$ and $B=\beta_S$. Since $|\Z_s|={p\choose s}$, the numerator of (\ref{eq:general-rigorous}) becomes $s+\log{p\choose s}\asymp s\log\left(\frac{ep}{s}\right)$, which is the well-known minimax rate of sparse linear regression \citep{donoho1994minimax,ye2010rate,raskutti2011minimax}. The principle (\ref{eq:general-rigorous}) also applies to a more general row and column sparsity structure in matrix denoising \citep{ma2015volume}.

\textit{Dictionary learning.} Consider the model $\X_Z(B)=BZ\in\mathbb{R}^{n\times d}$ for some $Z\{-1,0,1\}^{p\times d}$ and $Q\in\mathbb{R}^{n\times p}$. Each column of $Z$ is assumed to be sparse. Therefore, dictionary learning can be viewed as sparse linear regression without knowing the design matrix. With the relations $\tau=(p,s)$, $\T=\{(p,s)\in[n\wedge d]\times[n]:s\leq p\}$ and $\Z_{p,s}=\{Z\in\{-1,0,1\}^{p\times d}:\max_{j\in[d]}|\supp(Z_{*j})|\leq s\}$, we have $\ell(\Z_{p,s})+\log|\Z_{p,s}|\asymp np+ds\log\frac{ep}{s}$, which is the minimax rate of the problem \citep{klopp2017structured}.

The principle (\ref{eq:general}) or (\ref{eq:general-rigorous}) actually holds beyond the framework of structured linear models. We give an example of sparse principal component analysis (PCA). Consider i.i.d. observations $X_1,...,X_n\sim N(0,\Sigma)$, where $\Sigma=V\Lambda V^T+I_p$ belongs to the following space of covariance matrices
\begin{eqnarray*}
\mathcal{F}(s,p,r,\lambda) &=& \Bigg\{\Sigma=V\Lambda V^T+I_p: 0<\lambda\leq \lambda_r\leq...\leq\lambda_1\leq\kappa\lambda, \\
&&\qquad V\in O(p,r), |\rowsupp(V)|\leq s\Bigg\},
\end{eqnarray*}
where $\kappa$ is a fixed constant. The goal of sparse PCA is to estimate the subspace spanned by the leading $r$ eigenvectors $V$. Here, the notation $O(p,r)$ means the set of orthonormal matrices of size $p\times r$, $\rowsupp(V)$ is the set of nonzero rows of $V$, and $\Lambda$ is a diagonal matrix with entries $\lambda_1,...,\lambda_r$. It is clear that sparse PCA is a covariance model and does not belong to the class of structured linear models. 
Despite that, it has been proved in \cite{cai2013sparse} that the minimax rate of the problem is given by
\begin{equation}
\inf_{\wh{V}}\sup_{\Sigma\in\mathcal{F}(s,p,r,\lambda)}\mathbb{E}\fnorm{\wh{V}\wh{V}^T-VV^T}^2\asymp\frac{\lambda+1}{\lambda} \frac{r(s-r)+s\log\frac{ep}{s}}{n}. \label{eq:minimax-SPCA}
\end{equation}
The minimax rate (\ref{eq:minimax-SPCA}) can be understood as the product of $\frac{\lambda+1}{\lambda}$ and $\frac{r(s-r)+s\log\frac{ep}{s}}{n}$. The second term $\frac{r(s-r)+s\log\frac{ep}{s}}{n}$ is clearly a special case of (\ref{eq:general}). The first term $\frac{\lambda+1}{\lambda}$ can be understood as the modulus of continuity between the squared subspace distance used in (\ref{eq:minimax-SPCA}) and the intrinsic loss function of the problem (e.g. Kullback-Leibler), because the principle (\ref{eq:general}) generally holds for an intrinsic loss function. In addition to the sparse PCA problem, the minimax rate that exhibits the form of (\ref{eq:general}) or (\ref{eq:general-rigorous}) can also be found in sparse canonical correlation analysis (sparse CCA) \citep{gao2015minimax,gao2017sparse}.

\section{Community detection} \label{sec:community-detection}

\subsection{Problem settings}

The problem of community detection is to recover the clustering labels $\{z(i)\}_{i\in[n]}$ from the observed adjacency matrix $\{A_{ij}\}$ in the setting of SBM (\ref{eq:SBM-k}). It has wide applications in various scientific areas. 
Community detection has received growing interests in past several decades. 
Early contributions to this area focused on various cost functions to find graph clusters, in particular those based on graph cuts or modularity \citep{girvan2002community,newman2002random,newman2010networks}. 
Recent research has put more emphases on fundamental limits and provably efficient algorithms. 

In order for the clustering labels to be identifiable, we impose the following conditions in addition to (\ref{eq:SBM-k}),
\begin{equation}
\min_{1\leq u\leq k}B_{uu}\geq p,\quad \max_{1\leq u<v\leq k}B_{uv}\leq q.\label{eq:assortative}
\end{equation}
This is referred to as the assortative condition, which implies that it is more likely for two nodes in the same cluster to share an edge compared with the situation where they are from two different clusters. Relaxation of the condition (\ref{eq:assortative}) is possible, but will not be discussed in this survey. Given an estimator $\wh{z}$, we consider the following loss function
$$\ell(\wh{z},z)=\min_{\pi\in S_k}\frac{1}{n}\sum_{i=1}^n\indc{\wh{z}(i)\neq \pi\circ z(i)}.$$
The loss function measures the misclassification proportion of $\wh{z}$. Since permutations of labels correspond to the same clustering structure, it is necessary to take infimum over $S_k$ in the definition of $\ell(\wh{z},z)$.

In ground-breaking works by \cite{mossel2015reconstruction,mossel2013proof,massoulie2014community}, it is shown that the necessary and sufficient condition to find a $\wh{z}$ that is positively correlated with $z$ (i.e. $\ell(\wh{z},z)\leq\frac{1}{2}-\delta$) when $k=2$ is $\frac{n(p-q)^2}{2(p+q)}>1$. Moreover, the necessary and sufficient condition for weak consistency ($\ell(\wh{z},z)\rightarrow 0$) when $k=O(1)$ is $\frac{n(p-q)^2}{2(p+q)}\rightarrow\infty$ \citep{mossel2014consistency}. Optimal conditions for strong consistency ($\ell(\wh{z},z)=0$) were studied by \cite{mossel2014consistency,abbe2016exact}. 
When $k=2$, it is possible to construct a strongly consistent $\wh{z}$ if and only if $n(\sqrt{p}-\sqrt{q})^2>2\log n$, and extensions to more general SBM settings were investigated in \cite{abbe2015community}. We refer the readers to a thorough and comprehensive review by \cite{abbe2017community} for those modern developments.

Here we will concentrate on the minimax rates and algorithms that can achieve them. 
We favor the framework of statistical decision theory to derive minimax rates of the problem because the results automatically imply optimal thresholds in both weak and strong consistency. To be specific, the necessary and sufficient condition for weak consistency is that the minimax rate converges to zero, and the necessary and sufficient condition for strong consistency is that the minimax rate is smaller than $1/n$, because of the equivalence between $\ell(\wh{z},z)<1/n$ and $\ell(\wh{z},z)=0$. In addition, the minimax framework is very flexible and it allows us to naturally extend our results to more general degree corrected block models (DCBMs).

\subsection{Results for SBMs}
\label{sec:SBM}

We first formally define the parameter space that we will work with,
$$\Theta_k(p,q,\beta)=\left\{\theta=\{B_{z(i)z(j)}\}\in\Theta_k: n_u(z)\in\left[\frac{n}{\beta k},\frac{\beta n}{k}\right], B\text{ satisfies (\ref{eq:assortative})}\right\},$$
where the notation $n_u(z)$ stands for the size of the $u$th cluster, defined as $n_u(z)=\sum_{i=1}^n\indc{z(i)=u}$. We introduce a fundamental quantity that determines the signal-to-noise ratio of the community detection problem,
$$I=-2\log\left(\sqrt{pq}+\sqrt{(1-p)(1-q)}\right).$$
This is the R\'{e}nyi divergence of order $1/2$ between $\text{Bernoulli}(p)$ and $\text{Bernoulli}(q)$.
The next theorem gives the minimax rate for $\Theta_k(p,q,\beta)$ under the loss function $\ell(\wh{z},z)$.

\begin{thm}[Zhang and Zhou \citep{zhang2016minimax}]\label{thm:minimax-community-detection}
Assume $\frac{nI}{k\log k}\rightarrow\infty$, and then
\begin{equation}
\inf_{\wh{z}}\sup_{\Theta_k(p,q,\beta)}\mathbb{E}\ell(\wh{z},z)=\begin{cases}
\exp\left(-(1+o(1))\frac{nI}{2}\right), & k=2, \\
\exp\left(-(1+o(1))\frac{nI}{\beta k}\right), & k\geq 3,
\end{cases}\label{eq:minimax-community-detection}
\end{equation}
where $1+Ck/n\leq \beta <\sqrt{5/3}$ with some large constant $C>1$ and $0<q<p<(1-c_0)$  with some small constant $c_0\in(0,1)$. In addition, if $nI/k=O(1)$, then we have $\inf_{\wh{z}}\sup_{\Theta_k(p,q,\beta)}\mathbb{E}\ell(\wh{z},z)\asymp 1$.
\end{thm}

Theorem \ref{thm:minimax-community-detection} recovers some of the optimal thresholds for weak and strong consistency results in the literature. When $k=O(1)$, weak consistency is possible if and only if $\inf_{\wh{z}}\sup_{\Theta_k(p,q,\beta)}\mathbb{E}\ell(\wh{z},z)=o(1)$, which is equivalently the condition $nI\rightarrow\infty$ \citep{mossel2014consistency}. Similarly, strong consistency is possible if and only if $\frac{nI}{2}>\log n$ when $k=2$ and $\frac{nI}{\beta k}>\log n$ when $k$ is not growing too fast \citep{mossel2014consistency,abbe2016exact}. Between the weak and strong consistency regimes, the minimax misclassification proportion converges to zero with an exponential rate.

To understand why Theorem \ref{thm:minimax-community-detection} gives a minimax rate in an exponential form, we start with a simple argument that relates the minimax lower bound to a hypothesis testing problem. 
We only consider the case where $3\leq k=O(1)$ and $nI\rightarrow\infty$ are satisfied, and refer the readers to \cite{zhang2016minimax} for the more general argument. 
We choose a sequence $\delta=\delta_n$ that satisfies $\delta=o(1)$ and $\log\delta^{-1}=o(nI)$. Then, we choose a $z^*\in[k]^n$ such that $n_u(z^*)\in\left[\frac{n}{\beta k}+\frac{\delta n}{k},\frac{\beta n}{k}-\frac{\delta n}{k}\right]$ for any $u\in[k]$ and $n_1(z^*)=n_2(z^*)=\ceil{\frac{n}{\beta k}+\frac{\delta n}{k}}$. Recall the notation $\mathcal{C}_u(z^*)=\{i\in[n]:z^*(i)=u\}$. Then, we choose some $\wt{\mathcal{C}}_1\subset\mathcal{C}_1(z^*)$ and $\wt{\mathcal{C}}_1\subset\mathcal{C}_1(z^*)$ such that $|\wt{\mathcal{C}}_1|=|\wt{\mathcal{C}}_2|=\ceil{n_1(z^*)-\frac{\delta n}{k}}$. Define
$$T=\wt{\mathcal{C}}_1\cup \wt{\mathcal{C}}_2 \cup \left(\cup_{u=3}^k \mathcal{C}_u(z^*)\right)\quad\text{and}\quad\mathcal{Z}_T=\left\{z\in[k]^n:z(i)=z^*(i)\text{ for all }i\in T\right\}.$$
The set $\mathcal{Z}_T$ corresponds to a sub-problem that we only need to estimate the clustering labels $\{z(i)\}_{i\in T^c}$. Given any $z\in\mathcal{Z}_T$, the values of $\{z(i)\}_{i\in T}$ are known, and for each $i\in T^c$, there are only two possibilities that $z(i)=1$ or $z(i)=2$. The idea is that this sub-problem is simple enough to analyze but it still captures the hardness of the original community detection problem. Now, we define the subspace
$$\Theta_k^0(p,q,\beta)=\left\{\theta\in\{B_{z(i)z(j)}\}\in\Theta_k: z\in \mathcal{Z}_T, B_{uu}=p, B_{uv}=q,\text{ for all }1\leq u<v\leq k\right\}.$$
We have $\Theta_k^0(p,q,\beta)\subset\Theta_k(p,q,\beta)$ by the construction of $\mathcal{Z}_T$. This gives the lower bound
\begin{equation}
\inf_{\wh{z}}\sup_{\Theta_k(p,q,\beta)}\mathbb{E}\ell(\wh{z},z)\geq \inf_{\wh{z}}\sup_{\Theta_k^0(p,q,\beta)}\mathbb{E}\ell(\wh{z},z)=\inf_{\wh{z}}\sup_{z\in\mathcal{Z}_T}\frac{1}{n}\sum_{i=1}^n\mathbb{P}\{\wh{z}(i)\neq z(i)\}.\label{eq:reduction}
\end{equation}
The last inequality above holds because for any $z_1,z_2\in\mathcal{Z}_T$, we have $\frac{1}{n}\sum_{i=1}^n\indc{z_1(i)\neq z_2(i)}=O(\frac{\delta k}{n})$ so that $\ell(z_1,z_2)=\frac{1}{n}\sum_{i=1}^n\indc{z_1(i)\neq z_2(i)}$. Continuing from (\ref{eq:reduction}), we have
\begin{eqnarray}
\nonumber \inf_{\wh{z}}\sup_{z\in\mathcal{Z}_T}\frac{1}{n}\sum_{i=1}^n\mathbb{P}\{\wh{z}(i)\neq z(i)\} &\geq& \frac{|T^c|}{n}\inf_{\wh{z}}\sup_{z\in\mathcal{Z}_T}\frac{1}{|T^c|}\sum_{i\in T^c}\mathbb{P}\{\wh{z}(i)\neq z(i)\} \\
\label{eq:ratio} &\geq& \frac{|T^c|}{n}\frac{1}{|T^c|}\sum_{i\in T^c}\inf_{\wh{z}(i)}\text{ave}_{z\in\mathcal{Z}_T}\mathbb{P}\{\wh{z}(i)\neq z(i)\}.
\end{eqnarray}
Note that for each $i\in T^c$,
\begin{eqnarray}
\nonumber && \inf_{\wh{z}(i)}\text{ave}_{z\in\mathcal{Z}_T}\mathbb{P}\{\wh{z}(i)\neq z(i)\} \\
\label{eq:to-test} &\geq& \text{ave}_{z_{-i}}\inf_{\wh{z}(i)}\left(\frac{1}{2}\mathbb{P}_{(z_{-i},z(i)=1)}\left(\wh{z}(i)\neq 1\right)+\frac{1}{2}\mathbb{P}_{(z_{-i},z(i)=2)}\left(\wh{z}(i)\neq 2\right)\right).
\end{eqnarray}
Thus, it is sufficient to lower bound the testing error between each pair $\left(\mathbb{P}_{(z_{-i},z(i)=1)},\mathbb{P}_{(z_{-i},z(i)=2)}\right)$ by the desired minimax rate in (\ref{eq:minimax-community-detection}). Note that $|T^c|\gtrsim \frac{\delta n}{k}$ with a $\delta$ that satisfies $\log\delta^{-1}=o(nI)$. So the ratio $|T^c|/n$ in (\ref{eq:ratio}) can be absorbed into the $o(1)$ in the exponent of the minimax rate.

The above argument leading to (\ref{eq:to-test}) implies that we need to study the fundamental testing problem between the pair $\left(\mathbb{P}_{(z_{-i},z(i)=1)},\mathbb{P}_{(z_{-i},z(i)=2)}\right)$. That is, given the whole vector $z$ but its $i$th entry, we need to test whether $z(i)=1$ or $z(i)=2$. This simple vs simple testing problem can be equivalently written as
\begin{equation}
\label{eq:hypo}
\begin{aligned}
& H_1: X\sim \bigotimes_{i=1}^{m_1}\text{Bern}\left(p\right) \otimes \bigotimes_{i=m_1+1}^{m_1+m_2}\text{Bern}\left(q\right) \\
& \quad \quad \quad \mbox{vs.}\quad
H_2: X\sim \bigotimes_{i=1}^{m_1}\text{Bern}\left(q\right) \otimes \bigotimes_{i=m_1+1}^{m_1+m_2}\text{Bern}\left(p\right).	
\end{aligned}
\end{equation}
The optimal testing error of (\ref{eq:hypo}) is given by the following lemma.
\begin{lemma}[Gao, Ma, Zhang and Zhou \citep{gao2016community}]\label{lem:testing-error}
Suppose that as $m_1\rightarrow\infty$, $1<p/q=O(1)$, $p=o(1)$, $\left|m_1/m_2-1\right|=o(1)$ and $m_1I\rightarrow\infty$, we have
$$\inf_{\phi}\left(\mathbb{P}_{H_1}\phi+\mathbb{P}_{H_2}(1-\phi)\right)=\exp\left(-(1+o(1))m_1I\right).$$
\end{lemma}
Lemma \ref{lem:testing-error} is an extension of the classical Chernoff--Stein theory of hypothesis testing for constant $p$ and $q$ (see Chapter 11 of \cite{cover2012elements}). The error exponent $m_1I$ is a consequence of calculating the Chernoff information between the two hypotheses in (\ref{eq:hypo}). In the setting of (\ref{eq:to-test}), we have $m_1=(1+o(1))m_2=(1+o(1))\frac{nI}{\beta k}$, which implies the desired minimax lower bound for $k\geq 3$ in (\ref{eq:minimax-community-detection}). For $k=2$, we can slightly modify the result of Lemma \ref{lem:testing-error} with asymptotically different $m_1$ and $m_2$ but of the same order. In this case, one obtains $\exp\left(-(1+o(1))\frac{m_1+m_2}{2}I\right)$ as the optimal testing error, which explains why the minimax rate in (\ref{eq:minimax-community-detection}) for $k=2$ does not depend on $\beta$.

The testing problem between the pair $\left(\mathbb{P}_{(z_{-i},z(i)=1)},\mathbb{P}_{(z_{-i},z(i)=2)}\right)$ is also the key that leads to the minimax upper bound. Given the knowledge of $z$ but its $i$th entry, one can use the likelihood ratio test to recover $z(i)$ with the optimal error given by Lemma \ref{lem:testing-error}. Inspired by this fact, \citet{gao2017achieving} considered a two-stage procedure to achieve the minimax rate (\ref{eq:minimax-community-detection}). In the first stage, one uses a reasonable initial label estimator $\wh{z}^0$. This serves as an surrogate for the true $z$. Then in the second stage, one needs to solve the estimated hypothesis testing problem
\begin{equation}
\left(\mathbb{P}_{(z_{-i}=\wh{z}^0_{-i},z(i)=1)},\mathbb{P}_{(z_{-i}=\wh{z}^0_{-i},z(i)=2)},...,\mathbb{P}_{(z_{-i}=\wh{z}^0_{-i},z(i)=k)}\right).\label{eq:essential-testing}
\end{equation}
Since this is an upper bound procedure, we need to select from the $k$ possible hypotheses. The solution, derived by \cite{gao2017achieving}, is given by the formula
\begin{equation}
\wh{z}(i) = \argmax_{u\in[k]}\left(\sum_{\{j:\wh{z}^0(j)=u\}}A_{ij} - \wh{\rho}n_u(\wh{z}^0)\right). \label{eq:refinement}
\end{equation}
The number $\wh{\rho}$ is a data-driven tuning parameter that has an explicit formula given $\wh{z}^0$ (see \cite{gao2017achieving} for details). The formula (\ref{eq:refinement}) is intuitive. 
For the $i$th node, its clustering label is given by the one that has the most connections with the $i$th node, offset by the size of that cluster multiplied by $\wh{\rho}$. This one-step refinement procedure enjoys good theoretical properties. 
As was shown in \cite{gao2017achieving}, the minimax rate (\ref{eq:minimax-community-detection}) can be achieved given a reasonable initialization such as regularized spectral clustering. In practice, after updating all $i\in[n]$ according to (\ref{eq:refinement}), one can regard the current $\wh{z}$ as the new $\wh{z}^0$, and refine the estimator using (\ref{eq:refinement}) for a second round. From our experience, this will improve the performance, and usually less than ten steps of refinement is more than sufficient.

The ``refinement after initialization" method is a commonly used idea in community detection to achieve exponentially small misclassification proportion. Comparable results as \cite{gao2017achieving} are also obtained by \cite{yun2014accurate,yun2016optimal}. In addition to the likelihood-ratio-test type of refinement, the paper \cite{zhang2017theoretical} shows that a coordinate ascent variational algorithm also converges to the minimax rate (\ref{eq:minimax-community-detection}) given a good initialization.

\subsection{Results for DCBMs}

As was observed in \cite{bickel09,zhao2012consistency}, 
SBM is not a satisfactory model for many real data sets. 
An interesting generalization of SBM that captures degree heterogeneity was proposed by \cite{dasgupta2004spectral,karrer2011stochastic}, called degree corrected block model (DCBM). 
DCBM assumes that $A_{ij}\sim \text{Bernoulli}(d_id_jB_{z(i)z(j)})$. The extra sequence of parameters $(d_1,...,d_n)$ models individual sociability of network nodes. 
This extra flexibility is important in real-world network data analysis.

However, the extra nuisance parameters $(d_1,...,d_n)$ impose new challenges for community detection. 
There are not many papers that extend the results of SBM in \cite{mossel2015reconstruction,mossel2013proof,massoulie2014community,mossel2014consistency,abbe2016exact,abbe2015community} to DCBM. A few notable exceptions are \cite{zhao2012consistency,chen2015convexified,gulikers2015impossibility,gulikers2017spectral}.

On the other hand, the decision theoretic framework can be naturally extended from SBM to DCBM, and the results automatically imply optimal thresholds for both weak and strong consistency.

We first define the parameter space of DCBMs as
$$\Theta_k(p,q,\beta,d;\delta)=\left\{\theta=\{d_id_jB_{z(i)z(j)}\}: \{B_{z(i)z(j)}\}\in \Theta_k(p,q,\beta), \left|\frac{\sum_{i\in\mathcal{C}_u(z)}d_i}{n_u(z)}-1\right|\leq\delta\right\}.$$
Note that the space $\Theta_k(p,q,\beta,d;\delta)$ is defined for a given $d\in\mathbb{R}^n$. This allows us to characterize the minimax rate of community detection for each specific degree heterogeneity vector. The inequality $\left|\frac{\sum_{i\in\mathcal{C}_u(z)}d_i}{n_u(z)}-1\right|\leq\delta$ is a condition for $z$. This means that the average value of $\{d_i\}_{i\in\mathcal{C}_u(z)}$ in each cluster is roughly $1$, which implies approximate identifiability of $d,B,z$ in the model.

Before stating the minimax rate, we also introduce the quantity $J$, which is defined by the following equation,
\begin{equation}
\exp(-J) = \begin{cases}
\frac{1}{n}\sum_{i=1}^n\exp\left(-d_i\frac{nI}{2}\right), & k=2, \\
\frac{1}{n}\sum_{i=1}^n\exp\left(-d_i\frac{nI}{\beta k}\right), & k\geq 3.
\end{cases}\label{eq:J-def}
\end{equation}
When $d_i=1$ for all $i\in[n]$, $J$ is the exponent that appears in the minimax rate of SBM in (\ref{eq:minimax-community-detection}).

\begin{thm}[Gao, Ma, Zhang and Zhou \citep{gao2016community}]\label{thm:minimax-DCBM}
Assume $\min(J,\log n)/\log k\rightarrow\infty$, the sequence $\delta=\delta_n$ satisfies $\delta=o(1)$ and $\log\delta^{-1}=o(J)$, and $(d_1,...,d_n)$ satisfies Condition N in \cite{gao2016community}.
Then, we have
\begin{equation}
\inf_{\wh{z}}\sup_{\Theta_k(p,q,\beta,d;\delta)}\mathbb{E}\ell(\wh{z},z)=\exp(-(1+o(1))J),\label{eq:minimax-DCBM}
\end{equation}
where $1+Ck/n\leq \beta <\sqrt{5/3}$ and $1<p/q\leq C$  with some large constant $C>1$. 
\end{thm}

With slightly stronger conditions, Theorem \ref{thm:minimax-DCBM} generalizes the result of Theorem \ref{thm:minimax-community-detection}. By (\ref{eq:J-def}) and (\ref{eq:minimax-DCBM}), the minimax rate of community detection for DCMB is an average of $\exp(-d_i\frac{nI}{2})$ or $\exp(-d_i\frac{nI}{\beta k})$, depending on whether $k=2$ or $k\geq 3$. A node with a larger value of $d_i$ will be more likely clustered correctly.

Similar to SBM, the minimax rate of DCBM is also characterized by a fundamental testing problem. With the presence of degree heterogeneity, the corresponding testing problem is
\begin{equation}
\label{eq:hypo-DCBM}
\begin{aligned}
& H_1: X\sim \bigotimes_{i=1}^{m_1}\text{Bern}\left(d_0d_ip\right) \otimes \bigotimes_{i=m_1+1}^{m_1+m_2}\text{Bern}\left(d_0d_iq\right) \\
& \quad \quad \quad \mbox{vs.}\quad
H_2: X\sim \bigotimes_{i=1}^{m_1}\text{Bern}\left(d_0d_i q\right) \otimes \bigotimes_{i=m_1+1}^{m_1+m_2}\text{Bern}\left(d_0d_i p\right).	
\end{aligned}
\end{equation}
Here, we use $0$ as the index of node whose clustering label is to be estimated. The optimal testing error of (\ref{eq:hypo-DCBM}) is given by the following lemma.
\begin{lemma}[Gao, Ma, Zhang and Zhou \citep{gao2016community}]\label{lem:testing-error-DCBM}
Suppose that as $m_1\rightarrow \infty$, $1<p/q=O(1)$, $p\max_{0\leq i\leq m_1+m_2}d_i^2=o(1)$, $|m_1/m_2-1|=o(1)$ and $\left|\frac{1}{m_1}\sum_{i=1}^{m_1}d_1-1\right|\vee\left|\frac{1}{m_2}\sum_{i=m_1+1}^{m_2}d_i-1\right|=o(1)$. Then, whenever $d_0m_1I\rightarrow\infty$, we have
$$\inf_{\phi}\left(\mathbb{P}_{H_1}\phi+\mathbb{P}_{H_2}(1-\phi)\right)=\exp\left(-(1+o(1))d_0m_1I\right).$$
\end{lemma}
A comparison between Lemma \ref{lem:testing-error-DCBM} and Theorem \ref{thm:minimax-DCBM} reveals the principle that the minimax clustering error rate can be viewed as the average minimax testing error rate.

Finally, we remark that the minimax rate (\ref{eq:minimax-DCBM}) can be achieved by a similar ``refinement after initialization" procedure to that in Section \ref{sec:SBM}. 
Some slight modification is necessary for the method to be applicable in the DCBM setting, and we refer the readers to \cite{gao2016community} for more details.

\subsection{Initialization procedures} 

In this section, we briefly discuss consistent initialization strategies so that we can apply the refinement step (\ref{eq:refinement}) afterwards to achieve the minimax rate. 
The discussion will focus on the SBM setting. The goal is to construct an estimator $\wh{z}^0$ that satisfies $\ell(\wh{z}^0,z)=o_{\mathbb{P}}(1)$ under a minimal signal-to-noise ratio requirement. We focus our discussion on the case $k=O(1)$. Then, we need a $\wh{z}^0$ that is weakly consistent whenever $nI\rightarrow\infty$. The requirement for $\wh{z}^0$ when $k\rightarrow\infty$ was given in \cite{gao2017achieving}.

A very popular computationally efficient network clustering algorithm is spectral clustering \citep{shi2000normalized,mcsherry2001spectral,von2007tutorial,von2008consistency,rohe2011spectral,chaudhuri2012spectral,coja2010graph}. 
There are many variations of spectral clustering algorithms. 
In what follows, we present a version proposed by \cite{gao2016community} that avoids the assumption of eigengap. The algorithm consists of the following three steps:
\begin{enumerate}
\item Construct an estimator $\wh{\theta}$ of $\theta=\mathbb{E}A$.
\item Compute $\wt{\theta}=\argmin_{\rank(\theta)\leq k}\fnorm{\theta-\wh{\theta}}^2$.
\item Apply $k$-means algorithm on the rows of $\wt{\theta}$, and record the clustering result by $\wh{z}^0$.
\end{enumerate}

The three steps are highly modular and each one can be replaced by a different modification, which leads to different versions of spectral clustering algorithms \citep{zhou2018analysis}.
The vanilla spectral clustering algorithm either chooses the adjacency matrix or the normalized graph Laplacian as $\wh{\theta}$ in Step 1. Then, the $k$-means algorithm will be applied on the rows of $\wh{U}$ instead of those of $\wt{\theta}$ in Step 2, where $\wh{U}\in\mathbb{R}^{n\times k}$ is the matrix that consists of the $k$ leading eigenvectors. In comparison, the choice of $\wt{\theta}$ in Step 2 can be written as $\wt{\theta}=\wh{U}\wh{\Lambda}\wh{U}^T$, where $\wh{\Lambda}$ is a diagonal matrix that consists of the $k$ leading eigenvalues of $\wh{\theta}$. 
Our modified Step 1 and Step 2 make the algorithm consistent when $nI\rightarrow\infty$ without an eigengap assumption.

The choice of $\wh{\theta}$ in Step 1 is very important. Before discussing the requirement we need for $\wh{\theta}$, we need to understand the requirement for $\wt{\theta}$. According to a standard analysis of the $k$-means algorithm (see, for example, \cite{lei2015consistency,gao2016community}), a smaller $\mathbb{E}\fnorm{\wt{\theta}-\theta}^2$ leads to a smaller clustering error of $\wh{z}^0$. Therefore, the least-squares estimator (\ref{eq:sparse-least-squares}) will be the best option for $\wt{\theta}$ because it is minimax optimal (Theorem \ref{thm:minimax-sparse-graphon}). However, there is no known polynomial-time algorithm to compute (\ref{eq:sparse-least-squares}). On the other hand, the low-rank approximation in Step 2 can be computed efficiently through eigenvalue decomposition, and it enjoys the risk bound
$$\mathbb{E}\fnorm{\wt{\theta}-\theta}^2=O\left(\mathbb{E}\opnorm{\wh{\theta}-\theta}^2\right).$$
This means it is sufficient to find an $\wh{\theta}$ that achieves the minimal risk in terms of the squared operator norm loss. In other words, we seek an optimal sparse graphon estimator $\wh{\theta}$ with respect to the loss function $\opnorm{\cdot}^2$. The fundamental limit of the problem is given by the following theorem.
\begin{thm}[Gao, Lu and Zhou \citep{gao2015rate}]\label{thm:graphon-operator}
For $n^{-1}\leq \rho\leq 1$ and $k\geq 2$, we have
$$\inf_{\wh{\theta}}\sup_{\theta\in\Theta_k(\rho)}\mathbb{E}\opnorm{\wh{\theta}-\theta}^2\asymp \rho n.$$
\end{thm}
The minimax rate can be achieved by the estimator proposed by \cite{chin2015stochastic}. Define the trimming operator $T_{\tau}:A\mapsto T_{\tau}(A)$ by replacing the $i$th row and the $i$th column of $A$ with $0$ whenever $\sum_{j=1}^nA_{ij}\geq\tau$. Then, we set $\wh{\theta}=T_{\tau}(A)$ with $\tau=C\frac{1}{n-1}\sum_{i=1}^n\sum_{j=1}^nA_{ij}$ for some large constant $C>0$. This estimator can be shown to achieve the optimal rate given by Theorem \ref{thm:graphon-operator} \citep{chin2015stochastic,gao2017achieving}.
When $\rho \gtrsim \frac{\log n}{n}$ or the graph is dense, the native estimator $\wh{\theta}=A$ also achieves the minimax rate. This justifies the optimality of the results in \cite{lei2015consistency} in the dense regime.

With $\wh{\theta}$ described in the last paragraph, the three steps in the algorithm are fully specified. It can be shown that $\ell(\wh{z}^0,z)=o(1)$ with high probability as long as $nI\rightarrow\infty$ \citep{gao2017achieving,gao2016community}.

Another popular version of spectral clustering is to apply $k$-means on the leading eigenvectors of the normalized graph Laplacian $L=D^{-1/2}AD^{-1/2}$, where $D$ is a diagonal degree matrix. The advantage of using graph Laplacian in spectral clustering is discussed in \cite{von2008consistency,sarkar2015role}. When the graph is sparse, it is important to use regularized version of graph Laplacian \citep{amini2013pseudo,qin2013regularized,joseph2016impact}, defined as $L_{\tau}=D_{\tau}^{-1/2}A_{\tau}D_{\tau}^{-1/2}$, where $(A_{\tau})_{ij}=A_{ij}+\tau/n$ and $D_{\tau}$ is the degree matrix of $A_{\tau}$. The regularization parameter plays a similar role as the $\tau$ in $T_{\tau}(A)$. Performance of regularized spectral clustering is rigorously studied by \cite{le2015sparse,gao2017achieving,le2017concentration}.

For DCBM, it is necessary to apply a normalization for each row of $\wt{\theta}$ in Step 2. Instead of applying $k$-means directly on the rows of $\wt{\theta}$, it is applied on $\left\{\wt{\theta}_{1*}/\|\wt{\theta}_{1*}\|,...,\wt{\theta}_{n*}/\|\wt{\theta}_{n*}\|\right\}$. Then, the dependence on the nuisance parameter $d_i$ will be eliminated at each ratio $\wt{\theta}_{i*}/\|\wt{\theta}_{i*}\|$. The idea of normalization is proposed by \cite{jin2015fast} and is further developed by \cite{lei2015consistency,qin2013regularized,gao2016community}.

Besides spectral clustering algorithms, another popular class of methods is semi-definite programming (SDP) \citep{cai2015robust,chen2015convexified,amini2018semidefinite}. 
It has been shown that SDP can achieve strong consistency with the optimal threshold of signal-to-noise ratio \citep{hajek2016achieving,hajek2016achieving2}. Moreover, unlike spectral clustering, the error rate of SDP is exponential rather than polynomial \citep{fei2018exponential}.

\subsection{Some related problems}

The minimax rates of community detection for both SBM and DCBM are exponential. A fundamental principle for such discrete learning problems is the connection between minimax rates and optimal testing errors. In this section, we review several other problems in the literature that share this connection.

\textit{Crowdsourcing.} 
In many machine learning problems such as image classification and speech recognition, we
need a large amount of labeled data. Crowdsourcing provides an efficient while inexpensive
way to collect labels through online platforms such as Amazon Mechanical Turk \citep{turk2010url}.

Though massive in amount, the crowdsourced labels are usually fairly noisy. The low quality
is partially due to the lack of domain expertise from the workers and presence of spammers. Let $\{X_{ij}\}_{i\in[m],j\in[n]}$ be the matrix of labels given by the $i$th worker to the $j$th item. The classical Dawid and Skene model \citep{dawid1979maximum} characterizes the $i$th worker's ability by a confusion matrix
\begin{equation}
\pi_{gh}^{(i)}=\mathbb{P}(X_{ij}=h|y_j=g),\label{eq:dawid-skene}
\end{equation}
which satisfies the probabilistic constraint $\sum_{h=1}^k\pi_{gh}^{(i)}=1$. Here, $y_j$ stands for the label of the $j$th item, and it takes value in $[k]$. Given $y_j=g$, $X_{ij}$ is generated by a categorical distribution with parameter $\pi_{g*}^{(i)}=(\pi_{g1}^{(i)},...,\pi_{gk}^{(i)})$. The goal is to estimate the true labels $y=(y_1,...,y_n)$ using the observed noisy labels $\{X_{ij}\}$.

With the loss function $\ell(\wh{y},y)=\frac{1}{n}\sum_{j=1}^n\indc{\wh{y}_j\neq y_j}$, 
it is proved by \cite{gao2016exact} that under certain regularity conditions, the minimax rate of the problem is
\begin{equation}
\inf_{\wh{y}}\sup_{y\in[k]^n}\mathbb{E}\ell(\wh{y},y)=\exp\left(-(1+o(1))mI(\pi)\right),\label{eq:minimax-crowdsourcing}
\end{equation}
where
$$I(\pi)=-\max_{g\neq h}\min_{0\leq t\leq 1}\frac{1}{m}\sum_{i=1}^m\log\left(\sum_{l=1}^k\left(\pi_{gl}^{(i)}\right)^{1-t}\left(\pi_{hl}^{(i)}\right)^t\right)$$
is a quantity that characterizes the collective wisdom of a crowd.

The fact that (\ref{eq:minimax-crowdsourcing}) takes a similar form as (\ref{eq:minimax-community-detection}) is not a coincidence. The crowdsourcing problem is essentially a hypothesis testing problem. For each $j\in[n]$, one needs to select from the $k$ hypotheses $\{H_g\}_{g\in[k]}$, with the data generating process associated with $H_g$ given by (\ref{eq:dawid-skene}).

\textit{Variable selection.}
Consider the problem of variable selection in the Gaussian sequence model $X_j\sim N(\theta_j,\sigma^2)$ independently for $j=1,...,d$. The parameter space of interest is defined as
$$\Theta_d(s,a)=\left\{\theta\in\mathbb{R}^d: \theta_j=0\text{ if }z_j=0, \theta_j\geq a\text{ if }z_j=1, z\in\{0,1\}^d\text{ and }\sum_{j=1}^dz_j=s\right\}.$$
For any $\theta\in\Theta_d(s,a)$, there are exactly $s$ nonzero coordinates whose values are greater than or equal to $a$.

With the loss function $\ell(\wh{z},z)=\frac{1}{s}\sum_{j=1}^d\indc{\wh{z}_j\neq z_j}$, the minimax risk of variable selection derived by \cite{butucea2015variable} is
\begin{equation}
\inf_{\wh{z}}\sup_{\Theta_d(s,a)}\mathbb{E}\ell(\wh{z},z)=\Psi_+(d,s,a),\label{eq:minimax-variable-selection}
\end{equation}
where the quantity $\Psi_+(d,s,a)$ is given by
$$\Psi_+(d,s,a)=\left(\frac{d}{s}-1\right)\Phi\left(-\frac{a}{2\sigma}-\frac{\sigma}{a}\log\left(\frac{d}{s}-1\right)\right)+\Phi\left(-\frac{a}{2\sigma}+\frac{\sigma}{a}\log\left(\frac{d}{s}-1\right)\right).$$
The notation $\Phi(\cdot)$ is the cumulative distribution function of $N(0,1)$. Obviously, the optimal estimator that achieves the above minimax risk is the likelihood ratio test between $N(0,\sigma^2)$ and $N(a,\sigma^2)$ weighted by the knowledge of sparsity $s$.

Despite the connection between variable selection and hypothesis testing, the minimax risk (\ref{eq:minimax-variable-selection}) has two distinct features. First of all, the loss function is the number of wrong labels divided by $s$ instead of the overall dimension $d$. This is because the problem has an explicit sparsity constraint, which is not present in community detection or crowdsourcing. Second, given the Gaussian error, one can evaluate the minimax risk (\ref{eq:minimax-variable-selection}) exactly instead of just the asymptotic error exponent. The result (\ref{eq:minimax-variable-selection}) can be extended to more general settings and we refer the readers to \cite{butucea2015variable}.

\textit{Ranking.}
Consider $n$ objects with ranks $r(1),r(2),...,r(n)\in[n]$. We observe pairwise interaction data $\{X_{ij}\}_{1\leq i\neq j\leq n}$ that follow the generating process $X_{ij}=\mu_{r(i)r(j)}+W_{ij}$. The goal is to estimate the ranks $r=(r(1),...,r(n))$ from the data matrix $\{X_{ij}\}_{1\leq i\neq j\leq n}$. A natural loss function for the problem is $\ell_0(\wh{r},r)=\frac{1}{n}\sum_{i=1}^n\indc{\wh{r}(i)\neq r(i)}$. However, since ranks have a natural order, we can also measure the difference $|\wh{r}(i)-r(i)|$ in addition to the indicator whether or not  $\wh{r}(i)=r(i)$. This motivates a more general $\ell_q$ loss function $\ell_q(\wh{r},r)=\frac{1}{n}\sum_{i=1}|\wh{r}(i)-r(i)|^q$ for some $q\in[0,2]$ by adopting the convention that $0^0=0$. In particular, $\ell_1(\wh{r},r)$ is equivalent to the well known Kendall tau distance \citep{diaconis1977spearman} that is commonly used for a ranking problem.

Rather than discussing the general framework in \cite{gao2017phase}, we consider a special model with $X_{ij}\sim N(\beta(r(i)-r(j)),\sigma^2)$. Then, the minimax rate of the problem in \cite{gao2017phase} is given by
\begin{equation}
\inf_{\wh{r}}\sup_{r\in\mathcal{R}}\mathbb{E}\ell_q(\wh{r},r)\asymp\begin{cases}
\exp\left(-(1+o(1))\frac{n\beta^2}{4\sigma^2}\right), & \frac{n\beta^2}{4\sigma^2}>1, \\
\left[\left(\frac{n\beta^2}{4\sigma^2}\right)^{-1}\wedge n^2\right]^{q/2}, & \frac{n\beta^2}{4\sigma^2}\leq 1.
\end{cases}\label{eq:minimax-ranking}
\end{equation}
The set $\mathcal{R}$ is a general class of ranks that allow ties (approximate ranking). The detailed definition is referred to \cite{gao2017phase}. If $\mathcal{R}$ is replaced by the set of all permutations (exact ranking without tie), then the minimax rate (\ref{eq:minimax-ranking}) will still hold after $\frac{n\beta^2}{4\sigma^2}$ being replaced by $\frac{n\beta^2}{2\sigma^2}$ \citep{chen-minimax-ranking}.

The rate (\ref{eq:minimax-ranking}) exhibits an interesting phase transition phenomenon. When the signal-to-noise ratio $\frac{n\beta^2}{4\sigma^2}>1$, the minimax rate of ranking has an exponential form, much like the minimax rate of community detection in (\ref{eq:minimax-community-detection}). In contrast, when $\frac{n\beta^2}{4\sigma^2}\leq 1$, the minimax rate becomes a polynomial of the inverse signal-to-noise ratio.

When $\frac{n\beta^2}{4\sigma^2}>1$, the difficulty of ranking is determined by selecting among the following $n$ hypotheses
$$\left(\mathbb{P}_{r_{-i},r(i)=1},\mathbb{P}_{r_{-i},r(i)=2},...,\mathbb{P}_{r_{-i},r(i)=n}\right),$$
for each $i\in[n]$. Since the error of the above testing problem is dominated by the neighboring hypotheses, i.e., we need to decide whether $\mathbb{P}_{r_{-i},r(i)=j}$ or $\mathbb{P}_{r_{-i},r(i)=j+1}$ is more likely, the minimax ranking is then given by the exponential form in (\ref{eq:minimax-ranking}), where the exponent is essentially the Chernoff information between $\mathbb{P}_{r_{-i},r(i)=j}$ and $\mathbb{P}_{r_{-i},r(i)=j+1}$.


\section{Testing network structure} \label{sec:test}

\newcommand{\bbP}{\mathbb{P}}
\newcommand{\pnav}{p_{\mathrm{av}}}

\subsection{Likelihood ratio tests for \ER~model vs.~stochastic block model} \label{sec:LR-test}
Let $A = A^T\in \{0,1\}^{n\times n}$ be the adjacency matrix of an undirected graph with no self-loop.
For any probability $p\in [0,1]$, let $\mathcal{G}_1(n, p)$ denote the \ER~model where for all $i<j$, $A_{ij}\stackrel{iid}{\sim} \mbox{Bernoulli}(p)$.
For any $p\neq q\in [0,1]$, let $\mathcal{G}_2(n, p, q)$ denote the following ``mixture'' of SBMs. 
First, for $i=1,\dots, n$, let $z(i)-1\stackrel{iid}{\sim} \mbox{Bernoulli}(1/2)$. 
Conditioning on the realization of the $z(i)$'s, for all $i<j$, 
\begin{equation*}
A_{ij} = A_{ji} \stackrel{ind}{\sim}
\begin{cases}
\mbox{Bernoulli}(p), & \quad \mbox{if $z(i) = z(j)$},\\
\mbox{Bernoulli}(q), & \quad \mbox{if $z(i) \neq z(j)$}.
\end{cases}
\end{equation*}
We start with the simple testing problem of 
\begin{equation}
\label{eq:ervssbm}
H_0: A\sim \mathcal{G}_1\left(n, \frac{p+q}{2}\right)
~~~~\mbox{vs.}~~~~
H_1: A\sim \mathcal{G}_2(n, p, q).
\end{equation}
As in the previous section, we allow $p$ and $q$ to scale with $n$.
Under the present setting, both null and alternative hypotheses are simple, 
and so the Neyman--Pearson lemma shows that the most powerful test is the likelihood ratio test.
In what follows, we review the structure of the likelihood ratio statistics of the testing problem \eqref{eq:ervssbm} in two different regimes determined by whether the average node degree $\frac{n}{2}(p+q)$ remains bounded or grows to infinity as the graph size $n$ tends to infinity.
For convenience, denote the null distribution in \eqref{eq:ervssbm} by $\bbP_{0,n}$ and the alternative distribution $\bbP_{1,n}$.
We also define the following measure on the separation of the alternative distribution from the null 
\begin{equation}
	\label{eq:t-sep}
t = \sqrt{\frac{n(p-q)^2}{2(p+q)}}\,.
\end{equation}
The nontrivial cases are when $t$ is finite.

\paragraph{The regime of bounded degrees.}
In the asymptotic regime where 
\begin{equation}
	\label{eq:asymp-bdd-degree}
	np = a \quad\mbox{and} \quad 
	nq = b \quad\mbox{are constants as $n\to \infty$,}
\end{equation}
\citet{mossel2015reconstruction} focused on the assortative case where $p > q$ and showed that the testing problem \eqref{eq:ervssbm} has the following phase transition:
\begin{itemize}
\item When $t < 1$, $\bbP_{0,n}$ and $\bbP_{1,n}$ are asymptotically mutually contiguous, and so there is no consistent test for \eqref{eq:ervssbm};

\item When $t > 1$, $\bbP_{0,n}$ and $\bbP_{1,n}$ are asymptotically singular (orthogonal) and counting the number of cycles of length $\lfloor \log^{1/4}{n} \rfloor$ in the graph leads to a consistent test.
\end{itemize}
Their proof relies on the coupling of the local neighborhood of a vertex in the random graph with a Galton--Watson tree, which in turn depends crucially on the assumption \eqref{eq:asymp-bdd-degree} that the expected degrees of nodes remain bounded as the graph size grows.

In this asymptotic regime, when $t < 1$, the asymptotic distribution of the log-likelihood ratio can be characterized by that of the weighted sum of counts of $m$-cycles in the graph for $m \geq 3$. 
As far as asymptotic distribution is concerned, we can stop the summation at $m_n = \lfloor \log^{1/4} n \rfloor$. 
For each positive integer $j$, let $X_j = X_{n,j}$ be the counts of $j$-cycles in the graph of size $n$.
Following the lines of \cite{mossel2015reconstruction}, one can actually show that when $t< 1$ and \eqref{eq:asymp-bdd-degree} holds,
one achieves the asymptotic power of the likelihood ratio test by rejects for large values of 
\begin{equation}
\label{eq:lc-bdd}
L_c = \sum_{j=3}^{m_n} \qth{X_{n,j}\log(1+\delta_j) - \lambda_j \delta_j}
\end{equation}
with 
\begin{equation*}
\lambda_j = \frac{1}{2j}\left( \frac{a+b}{2} \right)^{j}
\qquad \mbox{and}\qquad
\delta_j = \left( \frac{a-b}{a+b} \right)^j.
\end{equation*}
To see this, one may replace Theorem 6 in \cite{mossel2015reconstruction} with Theorem 1 in \cite{Jan}.
Intuitively speaking, the counts of short cycles determine the likelihood ratio in the contiguous regime.

\paragraph{The regime of growing degrees.}
Now consider the following growing degree asymptotic regime where 
\begin{equation}
\label{eq:asymp-grow-degree}
np,\,nq\to\infty \quad \mbox{as}\quad n\to\infty.
\end{equation}
For simplicity, further assume that $p,q\to 0$ as $n\to \infty$, though all results in this part can be generalized to cases where $p$ and $q$ converge to constants in $(0,1)$.
Generalizing the contiguity arguments developed by \citet{Jan}, \citet{banerjee2018contiguity} established under \eqref{eq:asymp-grow-degree} the following phase transition phenomenon similar to that in the bounded degree case:
\begin{itemize}
\item When $t < 1$, $\bbP_{0,n}$ and $\bbP_{1,n}$ are asymptotically mutually contiguous, and so there is no consistent test for \eqref{eq:ervssbm};

\item When $t > 1$, $\bbP_{0,n}$ and $\bbP_{1,n}$ are asymptotically singular (orthogonal) and there is a consistent test.
\end{itemize}
Let $L_n = \frac{\mathrm{d} \bbP_{1,n}}{\mathrm{d} \bbP_{0,n}}$ be the likelihood ratio of \eqref{eq:ervssbm}.
\citet{banerjee2018contiguity} showed that when $t< 1$ and \eqref{eq:asymp-grow-degree} holds, the log-likelihood ratio $\log(L_n)$ satisfies 
\begin{equation}
\label{eq:llr-asymp}
\begin{aligned}
\log(L_n)  \stackrel{d}{\to} N\pth{-\frac{1}{2}\sigma(t)^2,\,\sigma(t)^2},
\qquad &\mbox{under $H_0$},\\
\log(L_n)  \stackrel{d}{\to} N\pth{\frac{1}{2}\sigma(t)^2,\,\sigma(t)^2},
\qquad &\mbox{under $H_1$},
\end{aligned}
\end{equation}
where 
\begin{equation*}
\sigma(t)^2 = \frac{1}{2}\pth{ -\log(1-t^2) - t^2 - \frac{t^4}{2} }.
\end{equation*}
Moreover, let $\pnav = \frac{1}{2}(p+q)$ and for any integer $i \geq 3$ define the \emph{signed cycles of length $i$} as
\begin{equation}
\label{eq:signedcycle}	
C_{n,i}(A) = \sum_{j_0, j_1,\dots, j_{i-1}} 
\qth{ \frac{A_{j_0 j_1}- \pnav} {\sqrt{n\pnav(1-\pnav)}}}\cdots 
\qth{ \frac{A_{j_{i-1} j_0}- \pnav} {\sqrt{n\pnav(1-\pnav)}}},
\end{equation}
where $j_0, j_1,\dots, j_{i-1}$ are all distinct and the summation is over all such $i$-tuples.
Unlike the actual counts of cycles used in the bounded degree regime, the signed cycles do not have a straightforward interpretation as graph statistics.
\citet{banerjee2018contiguity} further showed that when $t < 1$ and \eqref{eq:asymp-grow-degree} holds, the statistic
\begin{equation}
\label{eq:lc-stat}
L_{sc} = \sum_{i=3}^\infty \frac{2t^i C_{n,i}(A) - t^{2i}}{4i}
\end{equation}
has the same asymptotic distributions as those in \eqref{eq:llr-asymp} under both null and alternative.
In other words, a test that rejects $H_0$ for large values of $L_{sc}$ has the same asymptotic power as the likelihood ratio test which in turn is optimal by the Neyman--Pearson lemma.
Finally, when $t > 1$, with appropriate rejection regions, $L_{sc}$ leads to a consistent test.
Analogous to \eqref{eq:lc-bdd}, here the signed cycle statistics determine the likelihood ratio asymptotically within the contiguous regime.


%

\subsection{Tests with polynomial time complexity}

By our setting, the null hypothesis in \eqref{eq:ervssbm} is simple while the alternative averages over $2^n$ different possible configurations of the community assignment vector $z = (z(1),\dots, z(n))^T$.
Therefore, direct evaluation of the likelihood ratio test in either asymptotic regime is of exponential time complexity. 
It is therefore of great interest so see whether restricting one's attention to tests with polynomial time complexity would incur any penalty on statistical optimality \cite{berthet2013complexity, ma2015computational}. 
Interestingly, one can show that for the testing problem \eqref{eq:ervssbm}, there are polynomial time tests that are asymptotically as good as the likelihood ratio test.

\paragraph{The regime of bounded degrees.}
In view of \eqref{eq:lc-bdd}, in order to achieve the asymptotic powers of the likelihood ratio test, it suffices to count the numbers of $m$-cycles up to a slowly growing upper bound on $m$, say $m_n = \lfloor \log^{1/4} n \rfloor$.
Proposition 1 in \cite{mossel2015reconstruction} implied that this can be achieved within $\widetilde{O}(n(a+b)^{m_n})$ time complexity in expectation.

\paragraph{The regime of growing degrees.}
We divide the discussion into two different regimes.
In view of \eqref{eq:lc-stat}, it suffices to focus on estimating the signed cycles up to a slowly growing upper bound.
In what follows, we divide the discussion into two parts according to edge density.

First, assume that $np^2 \to \infty$. 
In this case, the average node degree grows at a faster rate than $\sqrt{n}$.
In this regime, \citet{banerjee2017optimal} showed that one can approximate the signed cycles by carefully designed linear spectral statistics of a rescaled adjacency matrix up to $m_n = \lfloor\min(\log^{1/2}(np^2), \log^{1/4} {n})\rfloor$.
In particular, define $A_{\mathrm{cen}} = (\frac{A_{ij} - \pnav}{\sqrt{n\pnav(1-\pnav)}}\mathbf{1}_{i\neq j} )$.
Moreover for any univariate function $g$ and any $n$-by-$n$ symmetric matrix $S$, let $\mathrm{Tr}(g(S)) = \sum_{i=1}^n g(\lambda_i)$ where the $\lambda_i$'s are the eigenvalues of $S$.
Furthermore, define $P_j(x) = 2S_j(x/2)$ where $S_j$ is the standard Chebyshev polynomial of degree $j$ given by $S_j(\cos\theta) = \cos(j\theta)$.
\citet{banerjee2017optimal} showed that when $np^2\to\infty$, a test that rejects for large values of 
\begin{equation}
	\label{eq:lc-chebyshev}
	L_a = \sum_{i=3}^{m_n} \frac{t^i}{2i}\mathrm{Tr}(P_i(A_\mathrm{cen}))
\end{equation}
achieves the asymptotic power of the likelihood ratio test within the contiguous regime.
Moreover, the mean and variance of $L_a$ admit explicit formulae, and so the computational cost of $L_a$ is $\widetilde{O}(n^3)$ as the most demanding step in its evaluation is computing the eigenvalues of an $n$-by-$n$ matrix.

Next, we consider the regime where $np\to\infty$, while $np^2$ remains bounded.
In this case, instead of working with the adjacency matrix directly, we may work with a scaled version of a weighted non-backtracking matrix proposed in \cite{fan2017well}.
For a graph with $n$ vertices, there are $n(n-1)$ distinct ordered pairs of $(i,j)$ with $i\neq j$.
Define a weighted non-backtracking matrix $B$ of size $(n^2-n)$-by-$(n^2-n)$ indexed by pairs of all such ordered pairs as
\begin{equation}
	\label{eq:weighted-nonbac}
B( (i,j), (i',j')  ) = 
\begin{cases}
	(A_\mathrm{cen})_{ij}, & \mbox{when $j=i'$ and $j'\neq i$,}\\
	0, & \mbox{otherwise}.
\end{cases}
\end{equation}
\citet{banma18} showed that as long as $np\to\infty$, a test based on some carefully constructed linear spectral statistic of $B$ achieves the asymptotic optimal power of the likelihood ratio test.
Regardless of the asymptotic condition, since the eigenvalues of $B$ can always be completed within $O(n^6)$ time complexity,
the time complexity of the test is bounded by $\widetilde{O}(n^6)$.



Finally, we mention that when $np\to\infty$ and $t>1$, \citet{montanari2016semidefinite} showed that SDP can be used to test \eqref{eq:ervssbm} consistently.


\subsection{Tests for more general settings}

When it comes to network data analysis in real world, degree heterogeneity is an indispensable feature for many social network data sets. This motivates us to consider a more general version of the hypothesis testing problem (\ref{eq:ervssbm}). We use $\mathcal{G}_k(n,p,q,\mathcal{D})$ to denote the following mixture of DCBMs. First, let $z(i)\stackrel{iid}{\sim}\text{Uniform}([k])$ for $i\in[n]$, and independently let $d_i\stackrel{iid}{\sim}\mathcal{D}$ for $i\in[n]$. Then, conditioning on $z(i)$'s and $d_i$'s, for all $i<j$,
\begin{equation*}
A_{ij} = A_{ji} \stackrel{ind}{\sim}
\begin{cases}
\mbox{Bernoulli}(d_id_jp), & \quad \mbox{if $z(i) = z(j)$},\\
\mbox{Bernoulli}(d_id_jq), & \quad \mbox{if $z(i) \neq z(j)$}.
\end{cases}
\end{equation*}
In order that the model parameters are identifiable, we impose the constraint
\begin{equation}
\mathbb{E}_{d\sim\mathcal{D}}(d^2)=1.\label{eq:2nd-moment-degree}
\end{equation}
When $p=q$ or $k=1$, the model is reduced to $A_{ij}|(d_i,d_j)\stackrel{ind}{\sim}\mbox{Bernoulli}(d_id_jp)$ for all $i<j$, which is recognized as as the configuration model \citep{van2016random} and is closely related to the Chung-Lu model \citep{chung2002average} of random graphs with expected degrees. The task we consider here is to test whether $p=q$ (or $k=1$) or $p\neq q$ (or $k>1$). The testing problem (\ref{eq:ervssbm}) can be viewed as a special case with $\mathcal{D}$ being a delta measure at $1$.

The key identity for the testing problem described above is revealed by the following lemma.
\begin{lemma}[Gao and Lafferty \citep{gao2017testing}]
Define the population edge (\edgeshape), vee (\veeshape), and triangle (\triangleshape) probabilities by $E=\mathbb{P}(A_{12}=1)$, $V=\mathbb{P}(A_{12}A_{13}=1)$, and $T=\mathbb{P}(A_{12}A_{13}A_{23})=1$. Then, under $A\sim \mathcal{G}_k(n,p,q,\mathcal{D})$ that satisfies (\ref{eq:2nd-moment-degree}), we have
\begin{equation}
T-\left(\frac{V}{E}\right)^3=\frac{(k-1)(p-q)^3}{k^3}.\label{eq:EZ}
\end{equation}
\end{lemma}
The relation (\ref{eq:EZ}) implies that $k=1$ or $p=q$ if and only if $T-(V/E)^3=0$. Intuitively speaking, when the network has more than one communities, its expected density of triangles deviates from the benchmark of a configuration model. When $T-(V/E)^3>0$, the network has an assortative clustering structure; such a network will induce more triangles compared with the configuration model. Conversely, the network will have a disassortative clustering structure if $T-(V/E)^3<0$, in which case there will be fewer triangles.

A remarkable feature of the equation (\ref{eq:EZ}) is its independence of the distribution $\mathcal{D}$ that characterizes the heterogeneity of the network nodes. Therefore, in order to test whether the null hypothesis is true or not, one does not need to estimate these nuisance parameters.

The asymptotic distribution of the empirical version of $T-(V/E)^3$ is given by the following theorem.
\begin{thm}[Gao and Lafferty \citep{gao2017testing}]\label{thm:GJ}
Consider the empirical versions of $E$, $V$, and $T$, defined as
\begin{eqnarray*}
\wh{E} &=& {n\choose 2}^{-1}\sum_{i<j}A_{ij}, \\
\wh{V} &=& {n\choose 3}^{-1}\sum_{i<j<l}\frac{A_{ij}A_{il}+A_{ij}A_{jl}+A_{il}A_{jl}}{3}, \\
\wh{T} &=& {n\choose 3}^{-1}\sum_{i<j<l}A_{ij}A_{il}A_{jl}.
\end{eqnarray*}
In addition to (\ref{eq:2nd-moment-degree}), assume $\mathbb{E}_{d\sim \mathcal{D}}(d^4)=O(1)$ and $n^{-1}\ll p\asymp q\ll n^{-2/3}$. Suppose
$$\delta=\lim_{n\rightarrow\infty}\frac{(k-1)(p-q)^3}{\sqrt{6}}\left(\frac{n}{k(p+(k-1)q)}\right)^{3/2}\in[0,\infty).$$
Then, we have
$$2\sqrt{{n\choose 3}}\left(\sqrt{\wh{T}}-(\wh{V}/\wh{E})^{3/2}\right)\leadsto N(\delta,1),$$
under the data generating process $A\sim \mathcal{G}_k(n,p,q,\mathcal{D})$.
\end{thm}
Under the null hypothesis, we have $k=1$, which implies $\delta=0$. This leads to the asymptotic distribution $N(0,1)$, and one can use the standard Gaussian quantile to determine the threshold of rejecting the null with a Type-1 error control. Under the alternative hypothesis, it is easy to see that $|\delta|\asymp \left(\frac{n(p-q)^2}{k^{4/3}(p+q)}\right)^{3/2}$. This implies a consistent test whenever $\frac{n(p-q)^2}{k^{4/3}(p+q)}\rightarrow\infty$, a condition that is slightly stronger than the optimal one discussed in Section \ref{sec:LR-test}, but applies to a much more general setting that even allows for a growing $k$.

One advantage of the above test is its applicability to real social network data because of both its simplicity and its invariance with respect to the distribution of the degree heterogeneity parameters.
\begin{figure}[!ht]
\begin{center}
\includegraphics[width=0.99\textwidth]{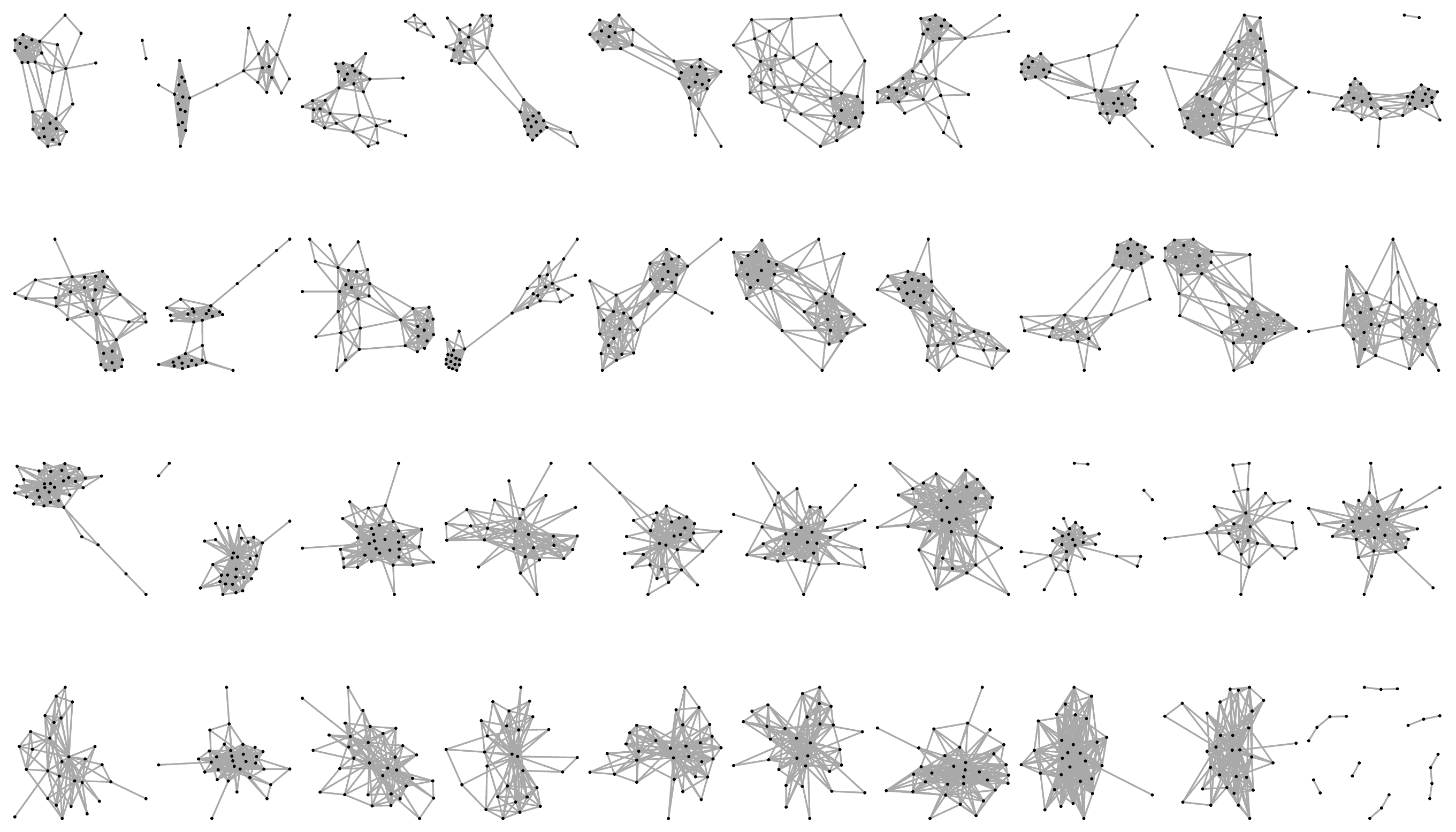}
\caption{Facebook neighborhood graphs, each with between 30 and 40 nodes,
  extracted from the Facebook 100 dataset. Top: 20 graphs with the smallest p-values. Bottom: 20 graphs
  with the largest. Community structure is readily apparent in
  the top graphs, and lacking in the bottom graphs.}\label{fig:fbgraphs}
\end{center}
\end{figure}
This provides practitioners a very useful tool to screen thousands of networks and only select those with potentially interesting community structure for further studies. Figure \ref{fig:fbgraphs} visualize real Facebook neighborhood graphs with small and large p-values.

Although the model $\mathcal{G}_k(n,p,q,\mathcal{D})$ is very flexible to derive a testing procedure that works really well in practice, further extensions are still possible by considering mixtures of degree corrected mixed membership models with possibly unbalanced community sizes. A test using more complicated subgraph counts including short paths and cycles is proposed by \cite{jin2018network}, and similar theoretical results as Theorem \ref{thm:GJ} are obtained.

\section{Discussion}\label{sec:disc}

There are a number of open problems in the research topics we have discussed.
For graphon estimation, the best rate achievable by polynomial-time algorithms has not been well understood. 
For community detection, even in the simple SBM setting, dependence of the tightest separation condition on the number of clusters $k$ is an interesting problem worth further investigation.
Furthermore, a challenging next step in network testing is to test a composite null of a $k$-community model against a composite alternative of models with more than $k$ communities.

The paper is focused on the notion of network sparsity introduced by \cite{bickel09,borgs2014p1,borgs2014p2}. Mathematically speaking, a network is sparse if $\max_{1\leq i<j\leq n}\theta_{ij}=o(1)$. However, this notion of sparsity contradicts the property of exchangeability \citep{lloyd2012random}, which is crucial for the inferential results to be able to generalize to the entire population \citep{mccullagh2002statistical}. Recently, two alternative notions of network sparsity have been developed in the literature. One of the proposals considers sparse networks induced by exchangeable random measures \citep{caron2017sparse,veitch2015class,borgs2016sparse}, and the other considers the notion of edge exchangeability \citep{crane2016edge,crane2018probabilistic}. Unlike the framework of \cite{bickel09,borgs2014p1,borgs2014p2}, these two alternative notions of sparsity allow well-defined sparse network models on the entire population, which implies a valid out-of-sample inference. However, rigorous and optimal statistical estimation and inference under these two frameworks are not well developed, except for only a few recent efforts \citep{todeschini2016exchangeable,herlau2016completely}. It is natural to ask whether the current state-of-the-art techniques of network analysis discussed in this paper can be modified or generalized to analyze sparse networks in these two alternative exchangeability frameworks. This question is of obvious significance and deserves extensive efforts of future research.

\begin{small}
\bibliographystyle{abbrvnat}
\bibliography{reference}
\end{small}


\end{document}